\documentclass[]{conm-p-l}
\usepackage{amsmath}
\usepackage{amssymb}
\usepackage{amscd}
\usepackage{amsthm}

\def\endproof{\relax
\ifmmode\expandafter\endproofmath
\else \unskip\nobreak\hfil\penalty50\hskip.75em\hbox{}\nobreak\hfil\bull
{\parfillskip=0pt \finalhyphendemerits=0
\bigbreak
}
\fi}
\def\endproofmath
$$
    {\eqno\bull
$$
\bigbreak
}
\def\bull{\vbox{\hrule\hbox{\vrule\kern3pt\vbox{\kern6pt}\kern3pt\vrule}\hrule}}

\newtheorem{theorem}{Theorem}[subsection]
\newtheorem{main}{Theorem}

\newtheorem{proposition}[theorem]{Proposition}
\newtheorem{lemma}[theorem]{Lemma}

\newtheorem{corollary}[theorem]{Corollary}

\newtheorem{D}[theorem]{Definition}
\newenvironment{defn}{\begin{D} \rm }{\end{D}}
\newtheorem{conjecture}[theorem]{Conjecture}
\newtheorem{R}[theorem]{Remark}
\newenvironment{remark}{\begin{R}\rm }{\end{R}}
\newcommand{%
\note}[1]{%
\marginpar{\scriptsize #1 }}

\def\Zee{\mathbb{Z}}

\def\Ar{\mathbb{R}}
\def\Cee{\mathbb{C}}
\def\Aff{\mathbb{A}}
\def\Pee{\mathbb{P}}

\def\Poin{\mathcal{P}}
\def\scrO{\mathcal{O}}
\def\Mod{\mathcal{M}}
\def\ov{\overline}
\def\spcheck{^{\vee}}
\def\frak{\mathfrak}
\def\Pic{%
\operatorname{Pic}
}
\def\Ker{%
\operatorname{Ker}
}
\def\Coker{%
\operatorname{Coker}
}
\def\Sym{%
\operatorname{Sym}
}
\def\Hom{%
\operatorname{Hom}
}

\def\Aut{%
\operatorname{Aut}
}
\def\Ext{%
\operatorname{Ext}
}
\def\Id{%
\operatorname{Id}
}

\def\Spec{%
\operatorname{Spec}
}
\def\ad{%
\operatorname{ad}
}
\def\Ad{%
\operatorname{Ad}
}
\def\End{%
\operatorname{End}
}
\def\Lie{%
\operatorname{Lie}
}
\def\Supp{%
\operatorname{Supp}
}

\begin{document}

\title[Minuscule Representations and Spectral Covers]{Minuscule Representations,
Invariant Polynomials, and Spectral Covers}

\author{Robert Friedman and John W. Morgan}
\thanks{The first author was partially supported by NSF grant
DMS-99-70437. 
The second author was partially supported by NSF grant DMS-97-04507.}

    \address{Department of Mathematics \\
    Columbia University \\
    New York, NY 10027}

    \email
        {rf@math.columbia.edu, jm@math.columbia.edu}


        \maketitle

\section*{Introduction}

Let $G$ be a simple and simply connected complex linear algebraic group, with Lie
algebra $\frak g$. Let $\rho\colon G \to \Aut V$ be an irreducible
finite-dimensional representation of $G$, and let $\rho_*\colon \frak g \to
\operatorname{End} V$ be the induced representation of $\frak g$. A  goal of
this paper is to study   $\rho_*$ and $\rho$, and in particular to give normal
forms for the action of
$\rho_*(X)$ and $\rho(g)$ for   regular elements $X$ of $\frak g$ or regular
elements $g\in G$. Of course, in the case of $\rho_*$, when
$X$ is semisimple, the action of $\rho_*(X)$  on $V$ can be diagonalized, and its
eigenvalues are given by evaluating the weights of $\rho$ with respect to a
Cartan subalgebra $\frak h$ of $\frak g$ on an element in  
$\frak h$ conjugate to
$X$. If on the other hand
$X$ is a principal nilpotent element of
$\frak g$, then $X$ can be completed to an $\frak{sl}_2$-triple $(X, h_0,X_-)$,
where
$h_0$ is regular and semisimple. In this case, if $\frak h$ is the Cartan
subalgebra containing
$h_0$, then the eigenvalues for   $\rho_*(h_0)$ together with their
multiplicities, which are given by evaluating the weights of $\rho$ with respect
to $\frak h$ on $h_0$, completely determine $V$ as an $\frak{sl}_2$-module and
hence determine the action of  $\rho_*(X)$  on $V$. A minor modification of these
ideas will then describe the action of $\rho_*(X)$ for every regular element $X$.

In this paper, we attempt to give a different algebraic model for the action
of $\rho_*(X)$, where $X$ is regular, and to glue these different models
together over the set of all regular elements.  Related methods
 also handle the case of
$\rho(g)$, where $g$ is a regular element of $G$. While we are only successful in
case
$\rho$ is minuscule, and partially successful in case
$\rho$ is quasiminuscule, it seems likely that the techniques of this paper can be
extended to give information about an arbitrary $\rho$. We believe that the
techniques and results of this paper are also of independent interest. For
example, they have been applied elsewhere to study sections of the adjoint
quotient morphism
\cite{Auto}.

To explain our results in more detail in the case of $\rho_*$, we begin by
recalling some results of Kostant on the adjoint quotient of $\frak g$. Let
$\frak h$ be the Lie algebra of a Cartan subgroup $H$ of $G$ and let $W$ be the
Weyl group of the pair
$(G,H)$. Kostant has shown \cite{Kostant} that the GIT quotient of $\frak g$ by
the adjoint action of
$G$ is isomorphic to $\frak h/W$. Moreover, he has constructed an explicit
cross-section $\ov\Sigma$ of the natural morphism $\frak g \to \frak h/W$, such
that the image of $\ov\Sigma$ is contained in the dense open subset $\frak
g_{\textrm{reg}}$ of regular elements of $\frak g$. 

 Given the Cartan subgroup $H$, there are the corresponding weights of $\rho$.
Fixing a choice of a set of positive roots for the pair $(G,H)$, we
shall let
$\mu$ be the
\textsl{lowest} weight of $\rho$.  Let $W_0$ be the stabilizer of $\mu$ in $W$.
Recall that  the representation $\rho$ is
\textsl{minuscule} if all of its weights are conjugate under $W$, so that the set
of  all weights of $\rho$ is identified with $W/W_0$. In this case, every weight
has multiplicity one. The representation $\rho$ is
\textsl{quasiminuscule} if all of its nonzero weights are conjugate under $W$.
In this case, every nonzero weight
has multiplicity one. 

Let $S=\Sym^*\frak h^*$ be the affine coordinate ring of $\frak h$, so that $S$ 
is a polynomial ring. We have the invariant
subrings $S^W$ and $S^{W_0}$. The inclusions $S^W \to S^{W_0} \to S$ correspond
to the finite morphisms $\frak h \to \frak h/W_0 \to \frak h/W$. By Chevalley's
theorem
\cite[p.\ 107, Thm.\ 3]{Bour},
$S^W$ and
$S^{W_0}$ are   polynomial algebras, i.e.\ $\frak h/W$ and $\frak h/W_0$ are
isomorphic to affine space. Hence
$S^{W_0}$ is flat over
$S^W$ and is therefore a free $S^W$-module of rank equal to $\#(W/W_0)$. Note that
$\mu$ is naturally an element of
$S^{W_0}$ and can be completed to a polynomial basis  of
$S^{W_0}$. The
morphism $\ov{\Sigma} \colon \frak h/W \to \frak g$ induces a morphism
$\rho_*\ov{\Sigma}\colon \frak h/W = \Spec  S^W  \to \operatorname{End} V$, and
is thus identified with an element of 
$$(\operatorname{End} V)\otimes _\Cee S^W\cong
\Hom_{S^W}(V\otimes _\Cee S^W,V\otimes _\Cee S^W),$$ 
which we  also denote
by $\rho_*\ov
\Sigma$. With this said, we have the following:

\begin{main}\label{main2} 
Let $\rho\colon G\to GL(V)$ be an irreducible finite-dimensional representation
with lowest weight
$\mu$. Then there is a nonzero linear map $f\colon V \to S^{W_0}$, unique up to a
nonzero scalar, with the following two properties: 
\begin{enumerate}
\item[\rm (i)] Let  $\hat f\colon
V\otimes_\Cee S^W\to S^{W_0}$ be the homomorphism defined by $f$ via extension of
scalars. Then the following diagram commutes:  
$$\begin{CD}
V\otimes_\Cee S^W @>{\hat f}>> S^{W_0} \\
@V{\rho_*\ov\Sigma}VV  @VV{\mu\cdot}V \\
V\otimes_\Cee S^W @>{\hat f}>> S^{W_0}.
\end{CD}
$$
\item[\rm (ii)] If $g\colon V \to S^{W_0}$ is any other linear map such that the
extension of $g$ to $\hat g\colon V\otimes_\Cee S^W\to S^{W_0}$ makes the above
diagram commute, then $g=sf$ for some $s\in S^{W_0}$.
\end{enumerate}
If $\rho$ is minuscule, then $\hat f$ is
an isomorphism.
If $\rho$ is quasiminuscule, then $\hat f$ is
surjective, and $\Ker \hat f =\Ker \rho_*\ov\Sigma$. 
\end{main}

In fact, we give an explicit description of the map $f$ in Section 3.

Using very different methods, and based on his results in \cite{Ginz}, Victor
Ginzburg has sent us the sketch of a proof that 
$\hat f$ is surjective for every irreducible representation $\rho$ (cf.\ also
\cite{Ginz2}).

There are two different ways we can interpret Theorem~\ref{main2}. In case
$\rho$ is minuscule, it describes, up to conjugation, the action of every regular 
element of
$\frak g$ on
$V$, or equivalently the action of $\rho_*\ov \Sigma(x)$ for every $x\in \frak
h/W$, in terms of the algebra structures on
$S^W$ and
$S^{W_0}$:

\begin{main}\label{action} Suppose that $\rho$ is minuscule. 
Let $x\in \frak h/W$ and let $\frak m_x\subseteq S^W$ be the maximal
ideal of  $x$, so that the scheme-theoretic fiber
over $x$ in $\frak h/W_0$
has coordinate ring equal to $S^{W_0}/\frak m_x\cdot S^{W_0}$. Then the map
$\hat f$ induces an isomorphism $\ov f_x\colon V\to S^{W_0}/\frak m_x\cdot
S^{W_0}$. Under this isomorphism, multiplication by $\mu$ on the
right-hand side becomes the action of  $\rho_*(\ov\Sigma(x))$. In particular, for
the ideal $\frak m_0\cdot S^{W_0}$ generated by the homogeneous $W$-invariant
functions of positive degree, there is an isomorphism
$\ov f_0\colon V\to S^{W_0}/\frak m_0\cdot S^{W_0}$ such that multiplication by
$\mu\in S^{W_0}$ on
$S^{W_0}/\frak m_0\cdot S^{W_0}$ corresponds to the action 
$\rho_*(X)$ of a principal nilpotent element $X$ on $V$.
\end{main}

Of course, the ring $S^{W_0}$ is just a polynomial ring and $\mu$ can be chosen
to be a coordinate function, so that the action of $\mu$ on $S^{W_0}$ is easy to
understand. However, the finite dimensional $\Cee$-algebras $S^{W_0}/\frak m_x
\cdot S^{W_0}$ and the action of $\mu$ on them are much more complicated. Thus for
example if $x=0$ the theorem tells us that the nilpotent element $\mu$ of 
$S^{W_0}/\frak m_0 \cdot S^{W_0}$ has Jordan blocks which are the same as
those given by the action of a principal nilpotent  element $X$ on $V$.

There is a partial analogue of Theorem~\ref{action} in the quasiminuscule,
non-minuscule case. In this case, the action of $\rho_*\ov \Sigma $ on $V\otimes
_\Cee S^W$ corresponds to the action of the matrix
$$\begin{pmatrix} 0&*\\0&\mu\end{pmatrix},$$
where we view $\mu$ as acting in block form on the free $S^W$-module $S^{W_0}$.
Here, for any given $x\in \frak h/W$, the off-diagonal term $*$ can in principle
be determined. We discuss this off-diagonal term in Section 5 for the adjoint
representation of a simply laced Lie algebra.

The second interpretation of the surjection $\hat f\colon V\otimes_\Cee S^W\to
S^{W_0}$ in the minuscule or quasiminuscule case is that it describes relations
between invariant polynomials for $W$ and $W_0$. In particular, an explicit
description  of
$f$ leads to a set of generators for $S^{W_0}$ as an $S^W$-module. We give
examples of such relations in the fourth section.

One basic ingredient in the proof of Theorem~\ref{main2} is a
result concerning the Kostant section $\ov\Sigma$ which (in a slightly different
form) was proved by Kostant in \cite{Kostant2}. In general, the image of
$\ov\Sigma$ is complicated. However, if we make the base change $\frak h\to \frak
h/W$, then, up to the adjoint action of a morphism
$\frak h\to G$,  the structure of $\ov \Sigma$ simplifies considerably:

\begin{main}\label{main1} Let $X$ be a principal nilpotent element of $\frak g$
such that $\frak h$ is contained in the unique Borel subalgebra containing $X$,
and let
$\Sigma$ be the composition of $\ov\Sigma$ with the natural projection $\frak
h\to \frak h/W$. Then there exists a morphism $\lambda\colon \frak h\to G$ such
that, for all $h\in \frak h$, 
$$\Ad(\lambda(h))(\Sigma (h)) = h+X.$$
\end{main}

In other words, after making the base change $\frak h\to \frak h/W$, we can
conjugate the section $\Sigma$ via the morphism $\lambda$ into the much simpler
morphism from $\frak h$ to $\frak g$ defined by $h\mapsto h+X$, where $X$ is a
fixed principal nilpotent in a Borel subalgebra normalized by $\frak h$. Using
the explicit description of the Kostant section, we are able to give a very
explicit construction of the morphism $\lambda$.

Steinberg has defined an analogous cross-section for the  quotient for the
adjoint action of $G$ on itself, which is a morphism $\ov \Phi \colon G\to H/W$
\cite{Steinberg}. We prove the following result concerning $\ov \Phi$, which is
analogous to Theorem~\ref{main1} but somewhat weaker:

\begin{main}\label{main4} Let $B$ be a Borel subgroup containing $H$ with
unipotent radical $U$, and let
$\Phi$ be the composition of $\ov\Phi$ with the natural projection $H\to H/W$.
Then there exists a morphism $\phi\colon H\to G$ and a morphism $u\colon
H\to U$ such that, for all
$h\in H$, 
$$\phi(h)\Phi (h)\phi(h)^{-1} = hu(h).$$
\end{main}

In contrast to the proof of Theorem~\ref{main1}, the argument does not use any
special properties of Steinberg's construction and would apply to any section of
the adjoint quotient morphism. In particular,  the morphism
$u$ is not given explicitly.  It is natural to conjecture that, for an appropriate
choice of
$\phi$ in the theorem, we can take $u$ to be a constant, in other words that the
pullback of the Steinberg section to
$H$ is conjugate via a morphism from $H$ to $G$ to the morphism $h\mapsto hu$ for
a certain principal unipotent element in a Borel subgroup normalized by $H$.  
Using  Theorem~\ref{main4}, the
proofs of Theorem~\ref{main2} and hence of Theorem~\ref{action} go over to the
case of the Steinberg section, and to the action of
$\rho(g)$ on $V$, for a regular element
$g$ in $G$ and a minuscule or quasiminuscule representation $\rho$. In particular,
if
$\mathcal{S}$ is the affine coordinate ring of $H$, then there is a morphism $\hat
g\colon V\otimes _\Cee\mathcal{S}^{W} \to \mathcal{S}^{W_0}$ which intertwines
the action of $\rho(\overline{\Phi})$ and multiplication by $\mu$, and $\hat g$ is
an isomorphism in case $\rho$ is  minuscule and surjective in case $\rho$ is
quasiminuscule. It is natural to conjecture that $\hat g$ is always surjective;
this does not seem to follow directly from the result of Ginzburg in the Lie
algebra case.

Our  strategy in the proof of Theorem~\ref{main2} is  to relate the problem
to a question concerning holomorphic $G$-bundles on a cuspidal curve. More
generally, let
$E$ be a Weierstrass cubic curve, i.e.\ a reduced irreducible curve of arithmetic
genus one. Then $E$ is either smooth, isomorphic to a nodal cubic curve in
$\Pee^2$, or isomorphic to a cuspidal cubic curve in $\Pee^2$. Let
$E_{\textrm{reg}}$ be the Zariski open subset of smooth points. Fixing a base
point
$p_0\in E_{\textrm{reg}}$ identifies
$E_{\textrm{reg}}$ with
$\Pic^0E$ and  endows $E_{\textrm{reg}}$ with the structure of a
one-dimensional connected, commutative algebraic group, which is  an elliptic
curve ($E$ smooth), or is isomorphic to the multiplicative group $\Cee^*$ ($E$
nodal)  or to the additive group 
$\Cee$ ($E$ cuspidal). We assume now that $E$ is singular and denote its
normalization by
$\widetilde E\cong
\Pee^1$. A
$G$-bundle on $E$ which pulls back to the trivial $G$-bundle $\widetilde E\times
G$ on $\widetilde E$ can be identified (after fixing a trivialization of the
pullback) with an element $g\in G$, if $E$ is nodal, or an element $X\in \frak
g$, if $E$ is cuspidal. There is a related result for families. In this way,
questions about $G$ or $\frak g$ can be translated into questions about
$G$-bundles $\xi$ on $E$. In case
$E$ is cuspidal, we can use the Kostant section to construct a universal
$G$-bundle $\Xi$ over $(\frak h/W)\times E$ which pulls back to the trivial
bundle over $(\frak h/W)\times \widetilde E$. A similar construction works in
the nodal case using the  Steinberg section. Although the construction of
$\Xi$ seems to depend on a particular choice of a section of the adjoint
quotient morphism, we show in \cite{Auto}, using many of the methods of this
paper, that all such sections are conjugate in an appropriate sense. If
$E$ is smooth, modulo some minor technical difficulties, we can use the parabolic
construction of
\cite{FMII}. Although the proof of Theorem~\ref{main2} only needs the case where
$E$ is singular, minor modifications also handle the case of a smooth $E$. This
case is interesting for the study of principal $G$-bundles on a smooth elliptic
curve $E$ as well as the vector bundles associated to them by a representation
$\rho$ of $G$. In fact, we were led to the results of this paper based
on our previous work on $G$-bundles over smooth elliptic curves $E$.

To prove Theorem~\ref{main2}, we shall use the theory of
spectral covers associated to semistable vector bundles over Weierstrass cubics. 
Given an irreducible representation
$\rho\colon G \to \Aut V$ and a $G$-bundle
$\xi$ on $E$, we can form the associated vector bundle $\xi\times_G V$. If $\xi$
pulls back to the trivial bundle on $\widetilde E$, then $\xi\times_G V$ is a
semistable vector bundle on $E$ of degree zero. It is well-known that, for every
Weierstrass cubic $E$, there is an equivalence between the abelian category of
semistable torsion free sheaves on $E$ of degree zero and the abelian category
of torsion sheaves on $E$ \cite{Mukai, FMW}. More precisely, given a semistable
torsion free sheaf
$V$ on $E$, there is a zero-dimensional subscheme $T_V$ in $E$, the
\textsl{spectral cover}, which is an effective Cartier divisor if $V$ is a
vector bundle, and a sheaf
$Q(V)$ supported on $T_V$ which can be used to recover $V$ from $Q(V)$ and  the
Poincar\'e sheaf on $E\times E$. The bundle $V$ is \textsl{regular} if $Q(V)$ is
a line bundle over $T_V$. More generally, given a family of semistable vector
bundles
$\mathcal{V}$ over $B\times E$, flat over $B$, the construction yields a relative
effective Cartier divisor
$T_{\mathcal{V}}$ in
$B\times E$ and a sheaf $Q(\mathcal{V})$ over $T_{\mathcal{V}}$, from which we
can recover $\mathcal{V}$. The bundle $\mathcal{V}$ restricts to a bundle which
is regular and semistable on every fiber if and only if $Q(\mathcal{V})$ is a
line bundle on
$T_{\mathcal{V}}$. However, we shall need to consider more general
situations, where $Q(\mathcal{V})$ is the pushforward of a line bundle on the
normalization of
$T_{\mathcal{V}}$.

In case $E$ is cuspidal, we use the Kostant section $\ov \Sigma$ to
construct a universal $G$-bundle $\Xi$ over $(\frak h/W)\times E$. Then we
let 
$\mathcal{V}$ be  the vector bundle over
$(\frak h/W)\times E$  associated to  $\Xi$ by the representation $\rho\colon G
\to \Aut V$.  Even though
$\Xi$ restricts to a regular
$G$-bundle
$\xi_x$ on every fiber
$\{x\}\times E$, the vector bundle $\mathcal{V}_x = \xi_x\times _GV$ is not in
general regular. The spectral cover $T_{\mathcal{V}}$ and the sheaf
$Q(\mathcal{V})$ are also quite complicated. For example, the irreducible
components of $T_{\mathcal{V}}$ are indexed by the Weyl orbits of the weights of
$\rho$, and each component has multiplicity equal to the multiplicity of the
orbit. However, the Weyl orbit of the highest weight (or equivalently of the
lowest weight) defines a reduced irreducible component $T'$ of $T_{\mathcal{V}}$.
Of course, $\rho$ is minuscule if and only if $T' = T_{\mathcal{V}}$. We
explicitly identify the normalization of $T'$ with a finite normal cover of
$\frak h/W$. The surjection $Q(\mathcal{V}) \to Q(\mathcal{V})|T'$ corresponds to
a surjection $\mathcal{V} \to \mathcal{V}'$, where $\mathcal{V}'$ is some coherent
sheaf on $(\frak h/W)\times E$. If $\rho$ is minuscule, then $\mathcal{V} =
\mathcal{V}'$. In the quasiminuscule case, we prove that $\mathcal{V}'$ is again
a vector bundle. It is natural to conjecture that this is true for an arbitrary
irreducible $\rho$. By analyzing how $\mathcal{V}'$ fails to be
regular in codimension one, we are able to establish Theorem~\ref{main2} but
without explicitly determining $\hat f$. To determine $\hat f$ requires a more
detailed analysis of the Kostant section which is given in Section 3. Similar
but less explicit constructions using the  Steinberg section work in case
$E$ is nodal.

The organization of this paper is as follows. In Section 1, we discuss bundles
over singular curves and recall the theory of spectral covers. We give  a
sufficient criterion for a vector bundle to be given by pushing forward a line
bundle on the normalization of the spectral cover.  In Section 2, we study the
spectral covers arising  from the vector bundles associated to a universal
$G$-bundle via an irreducible representation $\rho$, and in particular identify
the normalization of the irreducible component $T'$ associated to the Weyl orbit
of the highest weight. These results are then applied to the case of a minuscule
or quasiminuscule representation. In Section 3, we prove Theorem~\ref{main1}. We
can then easily prove Theorem~\ref{main2} by translating the vector bundle results
given in Section 2 into an algebraic form. The fourth section illustrates the
main results by studying various explicit examples. In Section 5, we study the
quasiminuscule, non-minuscule case in more detail in the simply laced case, where
the corresponding representation is the adjoint representation. It seems likely
that similar methods will also describe the non-simply laced case. We conclude
with some general conjectures concerning arbitrary representations.

It is a pleasure to thank Victor Ginzburg for explaining to us his work
\cite{Ginz, Ginz2} and for calling our attention to the reference
\cite{Kostant2}.

\section{Spectral covers of semistable vector bundles over Weierstrass cubics}

Throughout this paper, $E$ denotes a Weierstrass cubic, i.e.\ a reduced
irreducible curve of arithmetic genus one,   $E_{\textrm{reg}}$ is the set
of smooth points of $E$, and $p_0\in E_{\textrm{reg}}$ is a fixed point used
as the origin of the group law on $E_{\textrm{reg}}$. If $E$ is singular, then
its unique singular point is either a node or a cusp, and we shall refer to $E$
as either \textsl{nodal} or \textsl{cuspidal}. In this case, we denote the
normalization by
$\iota\colon
\widetilde E \to E$. The point
$p_0$ determines an isomorphism $E_{\textrm{reg}} \to \Pic^0E$. If $E$ is nodal,
$\Pic^0E$ is canonically $\Cee^*$, and if $E$ is cuspidal, the choice of a local
coordinate $z$ at $\iota^{-1}(E_{\textrm{sing}})$  determines an isomorphism of
$\Pic^0E$ with $\Cee$.  Let
$\Poin$ be the Poincar\'e bundle over $E\times E$. (For the meaning of this
sheaf and various properties of it in case $E$ is singular, we refer to the
proof of Lemma (0.3) in \cite{FMW}.) We shall use the symmetric form of
$\Poin$: $\Poin =
\scrO_{E\times E}(\Delta) \otimes
\pi_1^*\scrO_E(-p_0) \otimes \pi_2^*\scrO_E(-p_0)\otimes H^0(E;K_E)$, where
$\Delta \subseteq E\times E$ is the diagonal. In case
$E$ is smooth, the first factor is the parameter space $E=\Pic^0E$, and the second
factor is the original curve. Of course, the symmetry between the factors means
that we can also view the first factor as the curve and the second as
$\Pic^0E$. In the singular case, for the purposes of this paper, it will
suffice to consider the line bundle $\Poin|E_{\textrm{reg}}\times E$. Our goal in
this section is to describe basic results of
\cite{Mukai, FMW, Teo} which describe semistable bundles of degree zero on
$E$ in terms of torsion sheaves on
$E$. In \cite{FMW}, we mainly considered the case of bundles whose
restriction to every fiber was regular. Here we shall need a slight
generalization of those results. There are also analogous results relating
semistable torsion free sheaves of degree zero on
$E$ to torsion sheaves on $E$, but we shall not need them. 

\subsection{Bundles over singular curves}

Assume in this section that $E$ is singular and, in case $E$ is cuspidal, that
we have fixed once and for all a local coordinate $z$ on $\widetilde E$
centered at the inverse image of the singular point, and in particular have
determined an isomorphism $\Pic^0E\cong \Cee$. We begin by describing vector
bundles and
$G$-bundles over
$E$ which become trivial when pulled back to $\widetilde E$. More generally, we
shall consider parametrized versions. The arguments are standard, and left to
the reader. (See \cite{FMIII} for the case of a single curve.)

\begin{theorem}\label{triv} Suppose that $E$ is cuspidal and that we have fixed
the coordinate $z$ on $\widetilde E$ near the inverse image of the singular point
$x$. Let
$B$ be a scheme.   Suppose  that $\mathcal{V}$ is a vector bundle over $B\times E$
such that
$(\Id \times
\iota)^*\mathcal{V}$ restricts on every fiber to the trivial vector bundle of
rank $n$. In this case,
$\mathcal{W} = \pi_1{}_*(\Id \times \iota)^*\mathcal{V}$ is a rank $n$ vector
bundle on $B$, and $(\Id \times \iota)^*\mathcal{V} \cong \pi_1^*\mathcal{W}$.
Moreover, there  is  an   equivalence of categories between the 
category of triples $(\mathcal{V},
\mathcal{W},\Phi)$, where $\mathcal{V}$ is a rank $n$ vector bundle on
$B\times E$,
$\mathcal{W}$ is a rank $n$ vector bundle on $B$,  and
$\Phi$ is an isomorphism from $(\Id \times \iota)^*\mathcal{V}$ to
$\pi_1^*\mathcal{W}$, and the category of pairs $(\mathcal{W}, \varphi)$,
where 
$\mathcal{W}$ is a rank $n$ vector bundle on $B$ and $\varphi$ is an element
of $\operatorname{End} (\mathcal{W})$.
Similar statements hold in the nodal case, where $\varphi\in \Aut (\mathcal{W})$.
\qed
\end{theorem}

\begin{remark}\label{hartogsremark} (1) For $E$ cuspidal and $x =\iota^{-1}(E_{\rm
sing})$, the equivalence of categories is defined as follows: given
$(\mathcal{W}, \varphi)$, $\mathcal{V}$ is the subsheaf of $\pi_1^*\mathcal{W}$
consisting of all local sections $s$ such that, identically in $b$,
$$\frac{ds}{dz}(b,x) =
\varphi(b)s(b,x).$$
If $E$ is nodal, and $x,y\in \widetilde E$ are the two preimages of the singular
point on $E$, then 
$\mathcal{V}$ is the subsheaf of $\pi_1^*\mathcal{W}$
consisting of all local sections $s$ such that, for all $b$,
$$s(b,y) =
\varphi(b)s(b,x).$$
In particular, it follows that if $B$ is normal and $\mathcal{V}$ is a locally
free sheaf as in the statement of the theorem, then $\mathcal{V}$ has the Hartogs
property: if
$X$ is a subset of $B$ of codimension at least $2$ and $j\colon (B-X)\times E \to
B\times E$ is the inclusion, then the natural map
$\mathcal{V} \to j_*j^*\mathcal{V}$ is an isomorphism. In fact, a similar result
is true for every locally free sheaf on $B\times E$, as long as $B$ is normal.

\smallskip 
\noindent (2) The construction is compatible with base
change
$B'\to B$.

\smallskip 
\noindent (3) Suppose that $(\Id \times \iota)^*\mathcal{V}$ is  
trivialized. In this case, there is an equivalence of categories between the
category of triples
$(\mathcal{V}, V,\Phi)$, where
$\mathcal{V}$ is a vector bundle on $B\times E$, $V$ is a finite-dimensional
vector space,  and
$\Phi$ is an isomorphism from $(\Id \times \iota)^*\mathcal{V}$ to
$\scrO_{B\times
\widetilde E}\otimes _\Cee V$, and the category of pairs $(V,\varphi)$, where
$V$ is a finite-dimensional vector space and $\varphi\colon B \to 
\operatorname{End} V$ is a morphism.   A morphism $B\to \operatorname{End} V$ is 
the same  thing as a section of
$\scrO_B\otimes_\Cee \operatorname{End} V$. In case $B=\Spec R$, this is in turn 
equivalent to an element of $$\operatorname{End} V\otimes_\Cee R =
\operatorname{End} _R(V\otimes_\Cee R, V\otimes_\Cee R).$$ After choosing a basis
of $V$, this can further be identified with an
$n\times n$ matrix with coefficients in $R$. Similar remarks hold in the
nodal case. For example, in this case, $\varphi$ is a morphism
$ B \to \Aut V$, and in case $B=\Spec R$, such a morphism is
identified with an element of $\Aut_R(V\otimes _\Cee R)$.

\smallskip 
\noindent (4) The construction is natural under direct image by a finite
flat morphism. In other words, suppose that
$\nu\colon B
\to B'$ is a finite flat morphism. Thus $\nu_*\scrO_B$ is a locally free
$\scrO_{B'}$-module. If
$\mathcal{V}$ is a vector bundle on $B\times E$, then $(\nu\times
\Id)_*\mathcal{V}$ is a vector bundle on $B'\times E$. If in addition $(\Id
\times \iota)^*\mathcal{V}
\cong \pi_1^*\mathcal{W}$, then $(\Id \times \iota)^*(\nu\times
\Id)_*\mathcal{V}= \pi_1^*\nu_*\mathcal{W}$. In this case, the global section
$\varphi$ of $End (\mathcal{W})$ induces a global section of
$End(\nu_*\mathcal{W})$, and this corresponds to the bundle  $(\nu\times
\Id)_*\mathcal{V}$  on $B'\times E$.
\end{remark}

In case $B$ is a point and $\iota^*\mathcal{V}$ is trivialized, the rank $n$
bundle
$\mathcal{V}$ corresponds to an $n\times n$  matrix $A$ (which is arbitrary, in
the cuspidal case, and invertible, in the nodal case). For example, in the
cuspidal case, if
$N$ is the nilpotent matrix defined by $N(e_i) = e_{i-1}$, $i\geq 1$, and $N(e_1)
=0$, where the $e_i$ are the standard basis vectors, the corresponding vector
bundle is
$I_n$, in the notation of
\cite{FMW}. Here, for each $n>0$, $I_n$ is the 
unique indecomposable vector bundle  over $E$  of rank $n$ up to isomorphism which
has a filtration, all of whose successive quotients are $\scrO_E$. More
generally, if
$A = H+N$ is the Jordan decomposition of $A$, where $H$ is semisimple, $N$ is
nilpotent, and $[H,N] =0$, then the action of $H$ decomposes $V$ into a direct
sum of eigenspaces $\bigoplus _i V_i$, with distinct eigenvalues $c_i\in \Cee$,
and each $V_i$ is invariant under $N$. Suppose that the Jordan blocks of $N$
acting on $V_i$ have length
$n_{ij}$. Using the coordinate $z$ to identify $\Pic^0E$, we can identify
$c_i\in
\Cee$ with a line bundle on $E$ of degree zero, which we denote by $\lambda_i$.
Then
$$\mathcal{V}\cong \bigoplus _i\left(\bigoplus_jI_{n_{ij}}\right)\otimes
\lambda_i.$$
Of course, a similar statement holds for the nodal case. 

An important example in this paper will be the Poincar\'e bundle
$\Poin|E_{\textrm{reg}}\times E$. It is easy to check the following:

\begin{proposition}\label{cuspPoin}  Suppose that $E$ is cuspidal. The pullback
of $\Poin|E_{\rm reg}\times E$ to
$E_{\rm reg}\times \widetilde E$ is the trivial line bundle. Given a
coordinate $t$ on $E_{\rm reg}$ centered at $p_0$, under the
correspondence of Theorem~\ref{triv}, the Poincar\'e bundle 
$\Poin|E_{\rm reg}\times E$ corresponds to the morphism $\Cee \to \operatorname{End} \Cee$
given by multiplication, i.e.\ corresponds to the function $t$.\qed
\end{proposition}

A similar result handles the nodal case, where the coordinate $t$ on $E_{\rm
reg}\cong \Cee^*$ defines a isomorphism of algebraic groups from $E_{\rm reg}$ to
$\Aut\Cee$, again by multiplication.

We shall also use Theorem~\ref{triv} in the following situation. Let
$\rho\colon G \to \Aut V$ be a representation of $G$ and let $\rho_*\colon
\frak g\to \operatorname{End} V$ be the corresponding representation of the Lie algebra.
Given $X\in \frak g$, the element $\rho_*(X)$ defines a vector bundle over $E$.
If $X= X_s+X_n$ is the Jordan decomposition of $X$, then $\rho_*(X) = \rho_*(X_s)
+ \rho_*(X_n)$ is the Jordan decomposition of $\rho_*(X)$. Let  $\frak z_{\frak
g}(X)$ be the Lie algebra centralizer of $X$ in $\frak g$.  The element $X$ is
\textsl{regular} if and only if $\dim  \frak
z_{\frak g}(X)$ is equal to the rank of $\frak g$ if and only if $X_n$ is a
principal nilpotent element in
$\frak z_{\frak g}(X_s)$ Now suppose that we are given a morphism
$f\colon B\to
\frak g$,  for example  the   Kostant section  $\frak h/W \to
\frak g$ (which we shall describe later). Then
$\rho_*\circ f\colon B \to \operatorname{End} V$ is a morphism, and hence defines a vector
bundle $\mathcal{V}$ over $B\times E$ which becomes trivial over $B\times
\widetilde E$. Of course, a similar picture holds in the nodal case.

\subsection{Spectral covers}

\begin{defn} Let $V$ be a semistable bundle of degree zero on $E$. Define
$Q(V)$ as the cokernel of the natural global to local map $H^0(E; V\otimes
\scrO_E(p_0))\otimes_\Cee \scrO_E \to V\otimes_{\scrO_E} \scrO_E(p_0)$. Thus there
is an exact sequence
$$0 \to H^0(E; V\otimes \scrO_E(p_0))\otimes_\Cee \scrO_E \xrightarrow{A}
V\otimes_{\scrO_E} \scrO_E(p_0)
\to Q(V) \to 0.$$
We shall abbreviate $Q(V) = Q$ when $V$ is clear from the context, and call $Q$
the \textsl{Fourier-Mukai transform} of $V$.
\end{defn}

We have the following from \cite{FMW}:

\begin{lemma} Suppose that  $V$ is semistable of degree zero and rank $n$. Then
$h^0(E; V\otimes \scrO_E(p_0)) = n$,  $A$ is an injection from $\scrO_E^n$ to
$V\otimes \scrO_E(p_0)$, and $Q(V)$ is a torsion sheaf of length $n$.\qed
\end{lemma}

For example, if $L=\scrO_E(q-p_0)$ is a line bundle of degree zero, then
$Q(L) = \scrO_E/\frak m_q$. 

Let $T_V$ be the zero-dimensional subscheme of $E$
defined by the vanishing of the section $\det A$ of $\bigwedge ^n(V\otimes
\scrO_E(p_0))=\det V \otimes \scrO_E(np_0)$, in the above notation. By Cramer's
rule, $Q=Q(V)$ is   an
$\scrO_{T_V}$-module. Clearly $Q(V_1  \oplus V_2) = Q(V_1)\oplus Q(V_2)$. The
divisor $T_V$ is additive over  exact sequences: given $0\to V'\to V \to V''
\to 0$, where $V', V, V''$ are semistable of degree zero, we have $T_V =
T_{V'} + T_{V''}$ as divisors on $E$.

\begin{remark} (1) It follows from the defining exact sequence above that
$H^0(E;Q) \cong H^0(E; V\otimes \scrO_E(p_0)) \otimes H^1(E; \scrO_E)$ and
hence that there is a non-canonical isomorphism $H^0(E;Q) \cong H^0(E;
V\otimes \scrO_E(p_0))$.

\smallskip
\noindent (2) Mukai, in \cite{Mukai}, defines a torsion sheaf $Q'$ associated to
$V$ by: 
$$Q' = R^1\pi_1{}_*(\pi_2^*V\otimes \Poin^{-1}).$$
In fact, $Q' = Q\otimes
\scrO_E(p_0)$, which we can see as follows:   begin with the exact sequence
$$0 \to \Poin^{-1} \to  \pi_1^*\scrO_E(p_0) \otimes \pi_2^*\scrO_E(p_0) \to
\scrO_{\Delta}\otimes \pi_1^*\scrO_E(p_0) \otimes \pi_2^*\scrO_E(p_0) \to 0,$$
and tensor with $\pi_2^*V$. If we apply $R^i\pi_1{}_*$, we get an exact
sequence
$$0 \to H^0(E; V\otimes \scrO_E(p_0))\otimes \scrO_E(p_0) \to
V\otimes \scrO_E(2p_0)
\to R^1\pi_1{}_*(\pi_2^*V\otimes \Poin^{-1}) \to 0.$$
Hence $Q\otimes \scrO_E(p_0) \cong R^1\pi_1{}_*(\pi_2^*V\otimes \Poin^{-1})=Q'$.
\end{remark}

The \textsl{support} of $V$ is by definition the reduced support
of $T_V$, or equivalently of $Q$. If $\Supp V =\{q_1, \dots, q_a\}$, where
the $q_i$ are distinct, then $V=\bigoplus _iV(q_i)$, where each $V(q_i)$ is
supported at $q_i$. We shall only be concerned in this paper with bundles
$V$ such that the support of
$V$ is contained in
$E_{\textrm{reg}}$. We can see the structure
of $Q$ concretely in this case: If $\lambda$
is a line bundle of degree zero, let  $I_k(\lambda) = I_k\otimes \lambda$.
Given a point $q\in E_{\textrm{reg}}$, we let $\lambda_q =\scrO_E(q-p)$ be
the corresponding line bundle of degree zero. Then it is easy to check the
following:

\begin{lemma}\label{reg}  Suppose that $V\cong \bigoplus
_iI_{k_i}(\lambda_{q_i})$. Then
$T_V = \sum _ik_iq_i$ and 
$Q\cong \bigoplus _iQ_i$, where each $Q_i$ is supported at $q_i$ and in fact
$Q_i=\scrO_E/\frak m_{q_i}^{k_i}$. Thus $Q$ is a locally free
$\scrO_{T_V}$-module of rank one if and only if $q_i\neq q_j$ for
$i\neq j$, i.e.\ if and only if $V$ is regular. Finally, $T_V$ depends only
on the S-equivalence class of $V$. \qed
\end{lemma}

We have the following condition for a semistable vector bundle of degree
zero over a singular curve
$E$ to have support on $E_{\textrm{reg}}$:

\begin{lemma} Let $E$ be singular. Suppose that $V$ is a semistable vector bundle 
of rank $n$ and degree zero over $E$. Then the following are equivalent:
\begin{enumerate}
\item[\rm (i)] $\Supp V \subseteq E_{\rm reg}$;
\item[\rm (ii)] $V$ has a filtration whose successive quotients are line
bundles of degree zero;
\item[\rm (iii)] The pullback of $E$ to the normalization $\widetilde{E}$ is the
trivial rank
$n$ vector bundle.\qed
\end{enumerate}
\end{lemma}

\begin{lemma} Suppose that $E$ is cuspidal. Fixing a local coordinate at
the preimage in $\widetilde E$ of the singular point to identify
$E_{\rm reg}$ with $\Cee$,   let $t$ be the induced coordinate on $E_{\rm reg}$.
Let $V$ be a vector bundle over
$E$ whose pullback
$\widetilde{V}$ to $\widetilde{E}$ is trivial. Then, by Theorem~\ref{triv},  a
trivialization of 
$\widetilde{V}$ determines an element $A\in  \operatorname{End}
H^0(\widetilde{E};\widetilde{V})$. The conjugacy class of
$A$, viewed as an $n\times n$ matrix, is a complete invariant of the
isomorphism class of $V$. Finally, $T_V$ is the divisor in $E_{\rm reg}$ defined
by the characteristic polynomial $p_A(t)$ in $E_{\rm reg}$. A similar statement
holds if
$E$ is nodal.\qed
\end{lemma}

We return to the case of a general $E$ and $V$. To recover $V$ from $Q$, we
have the following result of Mukai
\cite{Mukai}:

\begin{proposition}\label{muk} With $Q=Q(V)$ defined as above, there is a
canonical isomorphism 
$$V\to \pi_2{}_*(\pi_1^*(Q\otimes
\scrO_E(p_0))\otimes \Poin).$$
\end{proposition}
\begin{proof} Let $H^0 = H^0(E; V\otimes \scrO_E(p_0))$. Using the exact
sequence
$$0\to H^0\otimes \scrO_E(p_0) \to V\otimes \scrO_E(2p_0) \to Q\otimes
\scrO_E(p_0) \to 0,$$
and the fact that $R^1\pi_2{}_*(H^0\otimes \pi_1^*\scrO_E(p_0))=0$, we see that it
suffices to find a canonical isomorphism 
$$\pi_2{}_*(\pi_1^*(V\otimes \scrO_E(2p_0))\otimes \Poin)\big
/\pi_2{}_*(\pi_1^*(H^0\otimes \scrO_E(p_0))\otimes \Poin) \to V .$$
To analyze the quotient, we begin with the exact sequence
$$0\to \scrO_{E\times E} \to \scrO_{E\times E}(\Delta) \to
\scrO_{\Delta}\otimes K_E^{-1} \to 0.$$
Using $\pi_1^*(V\otimes \scrO_E(2p_0))\otimes \Poin =
\pi_1^*(V\otimes \scrO_E(p_0))\otimes \pi_2^*\scrO_E(p_0)\otimes
\scrO_{E\times E}(\Delta)\otimes H^0(E;K_E)$, we get an exact sequence
$$0 \to \pi_1^*(V\otimes \scrO_E(p_0))\otimes \pi_2^*\scrO_E(-p_0) \to
\pi_1^*(V\otimes \scrO_E(2p_0))\otimes \Poin \to  i_*(V) \to
0,$$
where $i_*\colon \Delta \to E\times E$ is the inclusion. Applying $\pi_*$
then gives an exact sequence
$$0 \to H^0\otimes   \scrO_E(-p_0) \to \pi_2{}_*(\pi_1^*(V\otimes
\scrO_E(2p_0))\otimes \Poin) \to V \to 0.$$
As for the term $\pi_1^*(H^0\otimes \scrO_E(p_0))\otimes \Poin = H^0\otimes
\pi_2^*\scrO_E(-p_0)\otimes \scrO_{E\times E}(\Delta)$, the fact that
$\pi_2{}_*\scrO_{E\times E} \cong \pi_2{}_*\scrO_{E\times E}(\Delta)$ via the
natural inclusion shows that 
$$\pi_2{}_*(\pi_1^*(H^0\otimes
\scrO_E(p_0))\otimes \Poin) \cong H^0\otimes   \scrO_E(-p_0),$$
and it is easy
to check that this isomorphism is compatible with the inclusion of
$H^0\otimes   \scrO_E(-p_0)$ in $\pi_2{}_*(\pi_1^*(V\otimes
\scrO_E(2p_0))\otimes \Poin)$ given above. Thus we have identified 
$\pi_2{}_*(\pi_1^*(V\otimes \scrO_E(2p_0))\otimes \Poin)\big
/\pi_2{}_*(\pi_1^*(H^0\otimes \scrO_E(p_0))\otimes \Poin)$ with $V$.
\end{proof}

\begin{proposition}\label{inverseremark}  Given a torsion sheaf $Q$ on
$E$, if we set $$V(Q) = \pi_2{}_*(\pi_1^*(Q\otimes
\scrO_E(p_0))\otimes \Poin),$$ then $V(Q)$ is a semistable vector
bundle of degree zero, and there is a canonical isomorphism from the cokernel of
the natural map
$$H^0(E; V(Q)\otimes \scrO_E(p_0))\otimes \scrO_E \to V(Q)\otimes \scrO_E(p_0)$$
to $Q\otimes \scrO_E(p_0)$. 
\end{proposition} 
\begin{proof} Note
that
$V(Q) \otimes \scrO_E(p_0) \cong \pi_2{}_*(\pi_1^*Q\otimes \scrO_{E\times
E}(\Delta))$. 
Applying $\pi_2{}_*$ to the exact sequence
$$0 \to \pi_1^*Q\otimes \scrO_{E\times E}(\Delta - \pi_2^*p_0) \to 
\pi_1^*Q\otimes \scrO_{E\times E}(\Delta) \to \pi_1^*Q\otimes
\pi_2^*\scrO_E(p_0)|\pi_2^*p_0 \to 0$$
then gives the result.
\end{proof}

In other words,
$V(Q(V))\cong V$, and $Q(V(Q)) \cong Q\otimes \scrO_E(p_0)\cong Q$. Thus we have:

\begin{corollary} The functor $V\mapsto Q(V)$ is an
exact functor which defines an equivalence of categories from the abelian category
of semistable torsion free sheaves of degree zero over $E$ to the abelian category
of torsion  sheaves over
$E$. \qed
\end{corollary}

\begin{remark} Teodorescu \cite{Teo} gives an alternate method for recovering
$V$ from $Q$. It suffices to recover the extension class of the extension
$$0\to H^0\otimes \scrO_E \to V\otimes \scrO_E(p_0) \to Q \to 0,$$
where $H^0\otimes \scrO_E \cong H^0(E;Q)\otimes K_E$. But
$$\Ext^1(E; Q, H^0(E;Q)\otimes K_E) \cong \Ext^1(E; Q,K_E)\otimes H^0(E;Q)
\cong H^0(E;Q)^*\otimes H^0(E;Q),$$
by Serre duality, and the extension class is $\Id \in
H^0(E;Q)^*\otimes H^0(E;Q)$.
\end{remark}

\subsection{Families of vector bundles}

We turn now to the relative version of the construction. Throughout this section, 
$B$ denotes a scheme of finite type over $\Cee$ (or analytic space), and  
$\mathcal{V}$ is a rank
$n$ vector bundle over
$B\times E$. Given $b\in B$, let $\mathcal{V}_b$ be the bundle on $E$ induced by
the restriction of $\mathcal{V}$ to the slice $\{b\}\times E$. We shall always
assume that, for all $b\in B$, $\mathcal{V}_b$ is semistable of degree zero. We
can then define
$Q=Q(\mathcal{V})$ by the exact sequence
$$0\to \pi_1^*\pi_1{}_*(\mathcal{V}\otimes \pi_2^*\scrO_E(p_0))
\xrightarrow{A} \mathcal{V}\otimes \pi_2^*\scrO_E(p_0) \to Q \to 0,$$
and let $T_{\mathcal{V}}$ be the Cartier divisor in $B\times E$ defined by $\det
A$. We call $T_{\mathcal{V}}$ the \textsl{spectral cover} of $B$ defined by 
$\mathcal{V}$. Standard base change arguments show:

\begin{lemma}\label{relmuk} The morphism $T_{\mathcal{V}}\to B$ is finite and flat
of degree 
$n$. The coherent sheaf $Q$ on $B\times E$ is an
$\scrO_{T_{\mathcal{V}}}$-module, flat over $B$, and the rank of $Q$ viewed as a
vector bundle over $B$ is also $n$. Finally, the above construction is compatible
with base change: if $f\colon B'\to B$ is a morphism of schemes, then the
spectral cover associated to
$f^*\mathcal{V}$ is $f^*T_{\mathcal{V}}$ and the corresponding torsion sheaf in 
$B'\times E$ is $f^*Q$. \qed
\end{lemma}

For example, suppose that $L$ is a line bundle on $B\times E$ whose restriction
to every slice $\{b\}\times E$ has degree zero. Then $T_L$ is the section of
$B\times \Pic^0E= B\times E_{\rm reg}$ corresponding to $L$, and $Q(L)
=\scrO_{T_L}$.

Let $p_1=\pi_{1,2}\colon B\times E\times E \to B\times E$, let
$q=\pi_{1,3}\colon B\times E\times E \to B\times E$, and let
$r=\pi_{2,3}\colon B\times E\times E \to E\times E$. There is the relative
version of Proposition~\ref{muk}:

\begin{lemma} There is a natural isomorphism 
$$\mathcal{V} \cong
q_*(p_1^*(Q\otimes\pi_2^*\scrO_E(p_0))\otimes r^*\Poin).
\qed$$
\end{lemma}

Similarly, the analogue of  Proposition~\ref{inverseremark} holds: If 
$Q$ is a sheaf on $B\times E$, finite and flat over $B$, and $\mathcal{V} \cong
q_*(p_1^*(Q\otimes\pi_2^*\scrO_E(p_0))\otimes r^*\Poin)$, then 
$$Q(\mathcal{V})\cong  Q \otimes \pi_2^*\scrO_E(p_0).  
$$
In other words, up isomorphism, $Q(\mathcal{V})\cong  Q$.

We can form the restriction of
$r^*\Poin$ to $T_{\mathcal{V}}\times E$, which we shall denote 
$\Poin_{T_{\mathcal{V}}}$. Let $M = Q\otimes \pi_2^*\scrO_E(p_0)$, viewed as an
$\scrO_{T_{\mathcal{V}}}$-module, and let
$\nu\colon T_{\mathcal{V}}\to B$ be the projection. Then the lemma says that
$$\mathcal{V} \cong (\nu\times \Id)_*(\pi_1^*M\otimes 
\Poin_{T_{\mathcal{V}}}).$$

There are natural analogues of the above lemmas for flat families $\pi\colon Z\to
B$ of Weierstrass cubics with a section $\sigma$. Given a vector bundle
$\mathcal{V}$ over $Z$ whose restriction to every fiber of $\pi$ is
semistable of degree zero, we define $Q$ and $T_{\mathcal{V}}$ by the exact 
sequence
$$0\to \pi^*\pi _*(\mathcal{V}\otimes  \scrO_Z(\sigma))
\xrightarrow{A} \mathcal{V}\otimes  \scrO_Z(\sigma) \to Q \to 0,$$
and   $T_{\mathcal{V}}=\det A$,  and the obvious analogue of Lemma~\ref{relmuk}
holds.
Similarly, there is a relative Poincar\'e sheaf $\Poin =
\scrO_{Z\times_BZ}(\Delta)\otimes \pi_1^*\scrO_Z(-\sigma) \otimes
\pi_2^*\scrO_Z(-\sigma) \otimes p^*L$ over
$Z\times_BZ$, where $L$ is the line bundle on $B$ such that
$\omega_{Z/B}=\pi^*L$, where $\omega_{Z/B}$ is the relative dualizing sheaf
and $p\colon Z\times_B Z\to B$ is the projection, and a natural isomorphism
$$\mathcal{V} \cong \pi_2{}_*(\pi_1^*(Q\otimes \scrO_Z(\sigma))\otimes
\Poin).$$

We turn now to local properties of $Q$, or equivalently $M$ in the above
notation, viewed as an $\scrO_{T_{\mathcal{V}}}$-module. Applying Lemma~\ref{reg}, we have:

\begin{lemma} Suppose that $\mathcal{V}|\{b\}\times E$ has support in
$E_{\rm reg}$ for every $b\in B$. The
$\scrO_{T_{\mathcal{V}}}$-module
$Q$ is a line bundle if and only if the restriction of $\mathcal{V}$ to every
fiber is regular.\qed
\end{lemma}

As a  consequence of the lemma, suppose that $B$ is reduced and
equidimensional and that $\mathcal{V}|\{b\}\times E$ has support in
$E_{\rm reg}$ for every $b\in B$. Then $T$ is generically
reduced if and only if it is reduced, if and only if there is an open dense
subset
$U$ of
$B$ such that 
$\mathcal{V} |\{b\}\times E$ is isomorphic to a direct sum of $n$ distinct
line bundles of degree zero for all
$b\in U$, and in this case $Q$ is a torsion free rank one
$\scrO_{T_{\mathcal{V}}}$-module. (Here, $Q$ is \textsl{torsion free} if there 
are no  nonzero local sections of $Q$ supported on a subscheme of
$T_{\mathcal{V}}$ of strictly smaller dimension.) We shall need the following
more precise version:

\begin{proposition}\label{normalization} Suppose that
$B$ is normal and that $T_{\mathcal{V}}$ is reduced and irreducible. Let
$\hat{T}_{\mathcal{V}}$ be the normalization of $T_{\mathcal{V}}$, and let
$\hat\nu\colon \hat{T}_{\mathcal{V}} \to B$ be the composition
$\hat{T}_{\mathcal{V}}\to T_{\mathcal{V}}
\xrightarrow{\nu} B$. For each
$b\in B$, let
$\mathcal{V}_b$ be the bundle corresponding to $\mathcal{V}|\{b\}\times E$.
Suppose that there is a divisor
$D\subseteq B$ and an open dense subset $D'$ of $D$ such that
\begin{enumerate}
\item[\rm (i)] For $b\notin D$, $\mathcal{V}_b$ is regular and $T_{\mathcal{V}}$ 
is normal at all points of $(\nu)^{-1}(b)$.
\item[\rm (ii)] Given $b\in D'$ and $q\in \Supp \mathcal{V}_b$, exactly one
of the following holds:
\begin{enumerate}
\item[\rm (a)] 
$T_{\mathcal{V}}$ is normal at $(b,q)$ and $\mathcal{V}_b(q)$ is regular; 
\item[\rm (b)] 
$q\in E_{\rm reg}$, there exists a neighborhood   of $(b,q)$ in $T_{\mathcal{V}}$
of the form
$U_1\cup U_2$, where each $U_i$ is smooth, $(b,q)\in U_1\cap U_2$,
$U_1$ meets $U_2$ transversally, and
$\mathcal{V}_b(q) \cong \scrO_E(q-p_0)\oplus \scrO_E(q-p_0)$.
\end{enumerate}
\end{enumerate}
Then there exists a rank one reflexive sheaf $M$ on $\hat{T}_{\mathcal{V}}$ such
that 
$$\mathcal{V}
\cong (\hat\nu\times \Id)_*(\pi_1^*M\otimes \Poin_{\hat{T}_{\mathcal{V}}}),$$
where $\Poin_{\hat{T}_{\mathcal{V}}}$ is the
pullback of $\Poin_{T_{\mathcal{V}}}$ to $\hat{T}_{\mathcal{V}}\times E$. If
$\hat{T}_{\mathcal{V}}$ is smooth, then
$M$ is a line bundle.
\end{proposition}
\begin{proof} If $b\notin D$, then $T_{\mathcal{V}}$ is normal along 
$(\nu)^{-1}(b)$, so that $\hat{T}_{\mathcal{V}}=T_{\mathcal{V}}$ there, and $Q$
is locally free of rank one in a neighborhood of $(\nu)^{-1}(b)$. Now suppose
that $b\in D'$. If 
$T_{\mathcal{V}}$ is normal at $(b,q)$, then again  $Q$ is locally free of
rank one in a neighborhood of $(b,q)$. Suppose that
$T_{\mathcal{V}}$ is not normal at $(b,q)$. We first find locally a rank two
subbundle of $\mathcal{V}$ corresponding to the two branches of
$T_{\mathcal{V}}$. There is a small analytic Stein neighborhood
$U$ of
$b$ in $B$ such that $(\nu)^{-1}(U) = U'\cup U''$, where $(b,q)\in U'$, $U'$ is
connected, and
$U'\cap U'' =\emptyset$. By hypothesis, if $U$ is sufficiently small, we can
write 
$U'=U_1\cup U_2$, where
$U_1$ and
$U_2$ are smooth and meet transversally along some component $D_0$ of
$(\nu)^{-1}(D)$. Since $\nu$ is finite and flat and the local degree of $\nu$
at $(b,q)$ is
$2$, $\nu|U_i$ is an isomorphism of
$U_i$ onto $U$. Corresponding to $(\nu)^{-1}(U) = U'\cup U''$, we have
$\mathcal{V}|U\times E = \mathcal{V}' \oplus \mathcal{V}''$, where the rank of
$\mathcal{V}'$ is $2$. The restriction of $\Poin_{T_{\mathcal{V}}}$ to $U_i$ 
defines line bundles $L_i$ on $U_i\times E$ and hence on $U\times E$. Here $D_0$
is the preimage in $U_i$ of the set of points $u\in U$ where
$L_1|\{u\}\times E\cong L_2|\{u\}\times E$. Since
$(\mathcal{V}')\otimes L_1^{-1}|\{t\}\times E$ has a nonzero section
for each
$t\in U$, and has exactly one if $t\notin D$, it is easy to see that there is
a nonzero map $L_1\to \mathcal{V}'$. There is an effective divisor $Y$ on
$U\times E$ such that the map $L_1\to \mathcal{V}'$ extends to a map
$L_1\otimes \scrO_{U\times E}(Y)\to
\mathcal{V}'$ which does not vanish in codimension one. The quotient is
thus of the form $L_2'\otimes I_Z$, where $Z$ is a codimension two
subscheme of $U\times E$ (or empty).  Clearly
$\scrO_{U\times E}(Y)$ restricts to the trivial line bundle on every fiber
$\{u\}\times E$ and hence is pulled back from a line bundle on $U$. After
shrinking $U$, we may assume that this line bundle is trivial. There is a 
homomorphism from $L_2$ to $L_2'$ which is an isomorphism away from $D$. Thus
as before $L_2'\cong L_2\otimes
\scrO_{U\times E}(Y')$ and as before we may assume that the line bundle
$\scrO_{U\times E}(Y')$ is trivial. By semistability, the scheme
$Z$ is the preimage of a subset $Z'$ of
$U$, which is either empty or of codimension $2$. Thus in the complement  of
$Z'$,  there is an exact sequence
$$0 \to L_1\to \mathcal{V}' \to L_2 \to 0.$$
For $u\in U-Z'$, $\mathcal{V}'|\{u\}\times E$ is a direct sum of two line bundles
of degree zero. It follows that the extension
$$0 \to L_1|\{u\}\times E\to \mathcal{V}'|\{u\}\times E \to L_2 |\{u\}\times E\to
0$$
is split for every $u\in U-Z'$. For each $u\in U-Z'$, let  $U_u$ be
a Stein neighborhood of
$u$. Now 
$$\Ext^1(U_u\times E; L_2, L_1) = H^1(U_u\times E; (L_2)^{-1}\otimes L_1)) =
H^0(U_u; R^1\pi_2{}_*((L_2)^{-1}\otimes L_1)).$$
By the above discussion about $\mathcal{V}'|\{u\}\times E$, if $s$ is the
corresponding global section of $R^1\pi_2{}_*((L_2)^{-1}\otimes L_1)$, then the
image of
$s$ is zero in the fiber of $R^1\pi_2{}_*((L_2)^{-1}\otimes L_1)$ over $t$ for
every $t\in U_u$. To see that $s$ is actually zero, it suffices to prove that
$R^1\pi_2{}_*((L_2)^{-1}\otimes L_1)$ is supported along $D_0 = U_1\cap U_2$
and that its length is one there. Each line bundle $L_i$ defines a
morphism $g_i$ from $U$ to $\Pic^0E\cong E$, and the morphism $U_i\to
T_{\mathcal{V}}
\subseteq  U\times E$ is given by $(\Id, g_i)$. Let
$g=g_1-g_2$, so that $D_0 = g^{-1}(p_0)$. Then 
$(L_2)^{-1}\otimes L_1$ is the pullback to $U\times E$ via $(g\times \Id)^*$
of $\Poin$. The assumption that $U_1$ and $U_2$ meet transversally is the
assumption that
$g$ is smooth at the origin and hence that $D_0 = g^*(p_0)$ as reduced divisors.
By flat base change,
$$R^1\pi_2{}_*((L_2)^{-1}\otimes L_1)= R^1\pi_2{}_*(g\times
\Id)^*\Poin=g^*R^1\pi_1{}_*\Poin.$$
It is a standard result that $R^1\pi_1{}_*\Poin$ is supported at $p_0$ and
has length one there. It follows that $R^1\pi_2{}_*((L_2)^{-1}\otimes L_1)$ 
is supported along $D_0$ and has length   one there. Hence $s=0$. 

Thus
$\mathcal{V}'|U_u
\cong (L_1\oplus L_2)|U_u$. It follows that, at least over the complement of
the set $Z'$  of codimension at least two, 
$Q(\mathcal{V}')|(U-Z')\times E$, which is a direct summand of
$Q(\mathcal{V})|(U-Z')\times E$, is $Q(L_1) \oplus Q(L_2)|(U-Z')\times E \cong
\scrO_{U_1}\oplus
\scrO_{U_2}|(U-Z')\times E$. Hence, in the complement of the codimension two
subset $(\nu)^{-1}(Z')$ of 
$T_{\mathcal{V}}$,
$Q$ is a line bundle $M_0$ over the normalization
$\hat{T}_{\mathcal{V}}$ of $T_{\mathcal{V}}$. Since $\hat{T}_{\mathcal{V}}$ is
normal, $M_0$ extends uniquely to a rank one reflexive sheaf
$M$ on $\hat{T}_{\mathcal{V}}$. We have constructed   $M$  and an isomorphism
from $\mathcal{V}$ to $(\hat\nu\times \Id)_*(\pi_1^*M\otimes
\Poin_{\hat{T}_{\mathcal{V}}})$ over
$(B-Z')\times E$, where now $Z'$ is some global subscheme of $B$ of codimension
at least two. Since the codimension of $Z'$ in $B$ is at least two, if $j\colon
(B-Z')\times E \to B\times E$ is the inclusion, then it follows from (1) of
Remark~\ref{hartogsremark} that $\mathcal{V}\cong j_*j^*\mathcal{V}$. A similar
result holds for $(\hat\nu\times \Id)_*(\pi_1^*M\otimes
\Poin_{\hat{T}_{\mathcal{V}}})$, since $M$ is reflexive. Thus the isomorphism
from $\mathcal{V}$ to $(\hat\nu\times \Id)_*(\pi_1^*M\otimes
\Poin_{\hat{T}_{\mathcal{V}}})$ over
$(B-Z')\times E$
extends to give an isomorphism from 
$\mathcal{V}$ to $(\hat\nu\times \Id)_*(\pi_1^*M\otimes
\Poin_{\hat{T}_{\mathcal{V}}})$. The final statement follows since every rank one
reflexive sheaf on a smooth scheme is a line bundle.
\end{proof}

\subsection{Cohomology and extensions}

For applications to the study of quasiminuscule, non-minuscule representations,
we shall need a general result concerning extensions of a family of semistable
bundles by the trivial bundle. We begin with the case of a single  bundle over
$E$:

\begin{proposition} Let  $t$ be a local coordinate
centered at $p_0\in E$. Let $V$ be a torsion free
semistable sheaf on
$E$ of degree zero and let
$Q=Q(V)$. Then
\begin{enumerate}
\item[\rm (i)] $H^0(E; V) \cong \Ker \{\times t\colon Q \to Q\}$;
\item[\rm (ii)] $H^1(E; V) \cong \Coker \{\times t\colon Q \to Q\}=Q/tQ$.
\end{enumerate}
\end{proposition}
\begin{proof} To prove (i), let $\Cee_{p_0}$ be the torsion module corresponding
to
$\scrO_E$, i.e.\ 
$\Cee_{p_0} = \scrO_E/\frak m_{p_0} = \scrO_{E,p_0}/t\scrO_{E,p_0}$. Then, by
the equivalence of categories given by the Fourier-Mukai correspondence. 
$$H^0(E; V) = \Hom (\scrO_E, V) \cong \Hom(\Cee_{p_0}, Q) = \Ker \{\times t\colon
Q \to  Q\}.$$ To see (ii), the group $H^1(E; V)$ classifies the set of isomorphism
classes of extensions 
$0 \to V \to \mathcal{E} \to \scrO_E \to 0$,
and via the Fourier-Mukai correspondence this set is the same as the set 
of isomorphism classes of extensions
$0\to Q \to \tilde Q \to \Cee_{p_0} \to 0$,  which is classified by \
$\Ext^1_{\scrO_E}(\Cee_{p_0}, Q)$. But $\Ext^1_{\scrO_E}(\Cee_{p_0}, Q)\cong Q/tQ$
via the resolution 
$0\to \scrO_{E,p_0} \xrightarrow{\times t}  \scrO_{E,p_0} \to \Cee_{p_0} \to 0$.
\end{proof}

A minor extension of the argument to the relative case gives:

\begin{corollary} Let $B$ be a scheme, $\mathcal{V}\to B\times E$ a flat family of
torsion free semistable sheaves of degree zero and $Q=Q(\mathcal{V})$ the
corresponding sheaf on $B\times E$. Then, as sheaves on $B$, 
\begin{enumerate}
\item[\rm (i)] $R^0\pi_1{}_* \mathcal{V} \cong \Ker \{\times t\colon Q
\to
Q \}$;
\item[\rm (ii)] $R^1\pi_1{}_* \mathcal{V} \cong \Coker \{\times t\colon
Q \to Q\}=Q/tQ$. \qed
\end{enumerate}
\end{corollary}

As an application, we have:

\begin{corollary}\label{extcor} Let $B$ be a scheme, and let $\ov{\mathcal{V}}\to
B\times E$ a  family of torsion free semistable vector bundles of degree zero. Let
$\mathcal{W}$ be a vector bundle on
$B$. Then the relative extension sheaf
$Ext^1_{\pi_1}(\ov{\mathcal{V}},\pi_1^*\mathcal{W})
\cong R^1\pi_1{}_*(\ov{\mathcal{V}}\spcheck) \otimes \mathcal{W}$ is isomorphic to
$(Q\spcheck/tQ\spcheck)\otimes \mathcal{W}$, where
$Q\spcheck=Q(\ov{\mathcal{V}}\spcheck)$.
\qed
\end{corollary}

\section{Covers of the moduli space of $G$-bundles}

Throughout the remainder of this paper,  the rank of $G$ will be denoted by $r$,
$H$ is a maximal torus in $G$,
$\frak h=\Lie H$, and
$\Lambda \subseteq \frak h$ is the coroot lattice. The algebraic group
$E_{\textrm{reg}}\otimes_\Zee
\Lambda$,   is abstractly isomorphic to the $r$-fold
Cartesian product $(E_{\textrm{reg}})^r$. If
$E$ is cuspidal, then $E_{\textrm{reg}}\otimes_\Zee
\Lambda \cong \frak h$, and if $E$ is nodal then $E_{\textrm{reg}}\otimes_\Zee
\Lambda \cong H$. The Weyl group $W=W(G,H)$
acts on $\Lambda$. We let
$\Mod = (E_{\textrm{reg}}\otimes_\Zee \Lambda)/W$, and let $\widetilde \nu\colon E_{\textrm{reg}}\otimes_\Zee
\Lambda \to \Mod$ be the natural projection.
For
$E$ smooth,
$\Mod$ is the moduli space of semistable holomorphic $G$-bundles on $E$, modulo
S-equivalence. For $E$ nodal, $\Mod = H/W$ can be identified with the adjoint
quotient of $G$ by the conjugation action of $G$, where quotient is to be
taken in the GIT sense, i.e.\ again up to S-equivalence. By
\cite[\S6]{Steinberg}, $H/W$ is an affine space $\Aff^r$. Likewise,
if
$E$ is cuspidal, then $\Mod =\frak h/W$ is the (GIT) 
quotient of $\frak g$ via the adjoint action of $G$. By Chevalley's
theorem, $\frak h/W$ is again an affine space $\Aff^r$. (In case $E$ is nodal or
$E$ is cuspidal and $G$ is not of type $E_8$, the parabolic
construction gives a compactification of $\Mod$ \cite{FMIII}, which however will
not concern us here.)

\subsection{Definition of certain covers}

Fix a primitive weight $\mu$ for the pair $(G,H)$, i.e.\ a surjective homomorphism
$\Lambda
\to
\Zee$ (minor modifications handle the case where $\mu$ is not primitive). Let
$\Lambda_0 =\Ker\mu$ and let $W_0 =
\operatorname{Stab}_W(\mu)
\subseteq W$. Set $\hat{T}_0 = (E_{\textrm{reg}}\otimes_\Zee \Lambda)/W_0$. Thus
there is a finite morphism $\hat\nu\colon \hat{T}_0\to \Mod$. Since
$W_0$ acts on $\Lambda_0$, it also acts on $E_{\textrm{reg}}\otimes_\Zee
\Lambda_0$. We let $\Mod_0 = (E_{\textrm{reg}}\otimes \Lambda_0)/W_0$. The
homomorphism $\mu$ defines a $W_0$-invariant morphism
$E_{\textrm{reg}}\otimes_\Zee \Lambda\to E_{\textrm{reg}}$ and hence induces a
morphism $r\colon \hat{T}_0 \to E_{\textrm{reg}}$.

\begin{lemma} The space $\hat{T}_0$ is an \'etale fiber bundle over
$E_{\rm reg}$ with fibers isomorphic to $\Mod_0$.
\end{lemma}
\begin{proof} There exists a primitive $v\in \Lambda$ which is
$W_0$-invariant and whose image under $\mu$ is $d>0$. Define the surjection
$\Lambda\oplus \Zee \to \Zee$ by $(\lambda, k)\mapsto \mu(\lambda) - kd$.
Clearly, the kernel is exactly $\Lambda_0 \oplus \Zee$, included via the
standard inclusion on $\Lambda_0\oplus\{0\}$ and where $(0,1) \mapsto (v,1)$.
Tensoring with $E_{\textrm{reg}}$ gives a $W_0$-equivariant isomorphism
$E_{\textrm{reg}}\otimes_\Zee (\Lambda_0\oplus \Zee) =
(E_{\textrm{reg}}\otimes_\Zee
\Lambda_0)\times E_{\textrm{reg}} \to (E_{\textrm{reg}}\otimes_\Zee
\Lambda)\times _{E_{\textrm{reg}}}E_{\textrm{reg}}$, where in the second factor
the map $E_{\textrm{reg}}\to E_{\textrm{reg}}$ is multiplication by $d$. Taking
the quotient by
$W_0$ gives an isomorphism
$\hat{T}_0\times _{E_{\textrm{reg}}}E_{\textrm{reg}} \cong \Mod_0\times
E_{\textrm{reg}}$ as claimed.
\end{proof}

There is the following standard fact about groups generated by reflections:

\begin{lemma}\label{generated} The group $W_0$ is generated by the
reflections in $W$ which fix $\mu$, and hence by the reflections in the 
roots dual to the coroots lying in $\Ker\mu =\Lambda_0$. \qed
\end{lemma}

In the main application, $\mu$ will be a minuscule or quasiminuscule weight, and
in particular will be a fundamental weight $\varpi_\alpha$ for some simple root
$\alpha$ unless $G$ is of type $A_n$ and $\mu$ is a root. Clearly, if
$\mu = \varpi_\alpha$ for some  
$\alpha$, then 
$\Lambda_0$ is spanned by the simple coroots
$\beta$ for $\beta \neq
\alpha$. More generally, if $\Lambda_0$ is spanned by the coroots it
contains, then $\Lambda_0$ is the coroot lattice for a subroot system
$R_0$ of $R$ of rank $n-1$ and $W_0$ contains the Weyl group of $R_0$. In
fact, it follows from Lemma~\ref{generated} that $W_0$ is exactly the Weyl group
of $R_0$.  

Next we turn to the possible singularities of the moduli space and
their relation to the singularities of $\hat{T}_0$. Let $\Mod_{\textrm{sing}}$ be the
singular locus of $\Mod$, and let $\Mod_{\textrm{reg}} = \Mod
-\Mod_{\textrm{sing}}$. Note that $\Mod_{\textrm{reg}} = \Mod$ if $E$ is
singular,  or if $E$ is smooth and $G$ is of type $A_n$ or $C_n$, and in no
other cases. 

\begin{lemma} Suppose that $E$ is smooth. In the above notation,
$(\hat{T}_0)_{\rm{sing}} \subseteq \hat\nu^{-1}(\Mod_{\rm{sing}})$. 
\end{lemma}
\begin{proof} Let $K$ be the compact form of $G$. Viewing a point $e\in E 
\otimes_\Zee \Lambda$ as corresponding to a pair of elements $(x,y)\in T\times
T$, where $T$ is the maximal torus of $K$, let $K(e)$ be the centralizer of
$\{x,y\}$ in $K$. By
\cite[Theorem 6.12]{FMI}, the image of $e$ is a singular point of $\Mod$
if and only if $\pi_0(K(e)) \neq \{1\}$. Let $R_e$ be
the set of roots
$\alpha$ such that $e$ is in the kernel of the homomorphism
$E \otimes_\Zee \Lambda \to E$ defined by $\alpha$, and let  $W(R_e)$ be the
subgroup of
$W$ generated by reflections in the elements of $R_e$. A formal argument (cf.
\cite[Lemma 7.1.4]{BFM}) shows that there is an exact sequence
$$\{1\} \to  W(R_e) \to \operatorname{Stab}_W(e) \to \pi_0(K(e)) \to
\{1\}.$$
Then the image of
$e$ in
$\Mod$ is a singular point of $\Mod$ if and only if $\operatorname{Stab}_W(e)
\neq W(R_e)$.

Suppose that $e\in E\otimes_\Zee \Lambda$ lies over a smooth
point of $\Mod$. Then $\operatorname{Stab}_W(e) = W(R_e)$. Thus
$\operatorname{Stab}_{W_0}(e)$ is the stabilizer in $W(R_e)$ of $\mu$, and  
hence by Lemma~\ref{generated} is generated by the reflections in $W(R_e)$
which fix $\mu$. These generators fix $e\in E\otimes_\Zee \Lambda$ and act as
reflections on the tangent space to $E\otimes_\Zee \Lambda$ at $e$. Hence by
Chevalley's theorem
$\hat{T}_0$ is smooth in a neighborhood of the image of
$e$.
\end{proof}

\subsection{Existence of universal conformal bundles}

We now discuss the construction of universal $G$-bundles and universal
vector bundles over $\Mod \times E$. Let $E$ be a smooth elliptic curve.
Recall the following from
\cite{FMII}: for each $G$, there is 
\begin{itemize}
\item An integer $n>0$ and an element $c
\in Z(G)$ of order $n$;
\item A $G$-bundle $\Xi_0$ over $\Aff^{r+1}\times E$ and an action of
$\Cee^*$ on $\Xi_0$, lifting a linear action on $\Aff^{r+1}$,
\end{itemize}
with the following properties: 
\begin{itemize}
\item $\Xi_0|\{x\}\times E$ is semistable for all
$x\neq 0$;
\item the induced $\Cee^*$-invariant morphism $\Aff^{r+1}-\{0\} \to
\Mod$ induces an isomorphism from $\Aff^{r+1}-\{0\}/\Cee^*$ to $\Mod$.
\end{itemize}
Moreover, the primitive $n^{\rm th}$ roots of unity act trivially on
$\Aff^{r+1}$, the quotient $\Cee^*/(\Zee/n\Zee)$ acts effectively, and the
action of $\Zee/n\Zee$ on $\Xi_0$ is via $\langle c\rangle$. It follows
that, away from the singularities of $\Mod$, i.e.\ away from the non-free
$\Cee^*$-orbits, there is an induced $G/\langle c\rangle$ bundle.

In practice, we shall modify this construction slightly. Let $\hat G
= G\times _{\Zee/n\Zee}\Cee^*$, where $1\in \Zee/n\Zee$ acts via $c \in
G$ and as a primitive $n^{\rm th}$ root of unity $\zeta$ in $\Cee^*$. Given
$\zeta$, write $\zeta = \exp(2\pi\sqrt{-1}a/n)$, where $\gcd(a,n) =1$. Of
course, we could always fix $a=1$ if desired. There is a $\Cee^*$-linearization of
the trivial
$\Cee^*$-bundle over
$\Aff^{r+1}$, covering the given action on the base, as follows: $t(\lambda,
x) = (t^{-a}\lambda, t\cdot x)$. Taking the product bundle $\Cee^*\times
\Xi_0$ with its product linearization, we see that there is an induced $\hat
G$-bundle $\widehat\Xi_a$ over $\Mod_{\textrm{reg}}\times E$. 

\begin{lemma} For each $x\in \Mod_{\rm reg}$, the $\hat G$-bundle
$\widehat\Xi_a|\{x\}\times E$ has a canonical lift to a $G$-bundle $\xi_x$,
which is the regular $G$-bundle corresponding to the point $x$.
\end{lemma}
\begin{proof} Begin with the central extension
$$\{1\} \to G \to \hat G \to \Cee^* \to \{1\}.$$
By construction, the bundle $\hat \xi_x =\widehat\Xi_a|\{x\}\times E$ lifts to
a
$G$-bundle. Let $\det \hat \xi_x$ be the $\Cee^*$-bundle induced from the
quotient map $\hat G \to \Cee^*$. The inclusion
$\Cee^*\subseteq
\hat G$ defines a surjective homomorphism $\Aut \hat \xi_x \to \Aut \det\hat
\xi_x =\Cee^*$, and hence the lift of $\hat \xi_x$ is unique.
\end{proof}

For
$E$ nodal or cuspidal, we could use slight variations of this construction
(with special care in case
$E$ is cuspidal and $G=E_8$). In fact, the moduli space $\Mod$ is
isomorphic to an affine subspace $\Aff^r$ in $\Aff^{r+1}-\{0\}$ as we
show in \cite{FMIII}, and hence there is a universal
$G$-bundle $\Xi$ over $\Mod\times E$. On the other hand, the Steinberg and Kostant
sections already give us universal
$G$-bundles as described in \cite{FMIII}. As we shall show in \cite{Auto}, all
of these sections are conjugate in the appropriate sense, so that the choice
does not really matter. However, in the next section, for the cuspidal case, we
shall use the explicit form of the Kostant section.

Now suppose that $\rho\colon G \to GL(V)$ is an irreducible representation.
Then $\rho(c) = \zeta'\Id$, where $\zeta'$ is a root of unity whose order
divides $n$. Thus there is $b\in \Zee$ such that $\zeta^b =\zeta'$. The map
$\chi_b(t) = t^b\Id$ is a homomorphism from $\Cee^*$ to $GL(V)$, and
$\rho\otimes \chi_b$ factors to give a representation $\rho_b\colon \hat G
\to GL(V)$. There is then the associated bundle $\mathcal{V}_{a,b} =
\widehat\Xi_a\times _{\hat G}V$ over
$\Mod_{\textrm{reg}}\times E$, where $\hat G$ acts on
$V$ via $\rho_b$. It is easy to see that the isomorphism type of the vector bundle
$\mathcal{V}_{a,b}$ only depends on the product $ab$. For each $x\in
\Mod_{\textrm{reg}}$,$\mathcal{V}_{a,b}|\{x\}
\times E = \xi_x\times_GV  =V_x$, and hence does not depend on $a,b$. Since
the image of $\rho$ is contained in $SL(V)$, the bundle $V_x$ has trivial
determinant. Thus $\det \mathcal{V}_{a,b}$ is pulled back from the factor
$\Mod_{\textrm{reg}}$.  We can take the associated spectral cover
$T_{\mathcal{V}_{a,b}}$ in the sense of \S 1, and it is also independent of $a,b$.
Finally, it is not hard to see that $Q(\mathcal{V}_{a,b})$ is independent of the
choice of $(a,b)$ up to twisting by a line bundle on $T_{\mathcal{V}_{a,b}}$. Fix
some pair $(a,b)$, and set $\mathcal{V}=\mathcal{V}_{a,b}$.

The fiber of $T_{\mathcal{V}}$ over $x$ only depends on the S-equivalence class of
$\xi_x$. Thus we have:

\begin{lemma} Let $\mu_1, \dots, \mu_n$ be
the   weights of
$V$, counted with multiplicity. Suppose that $\xi_x$ is S-equivalent to the
split bundle corresponding to the point $e\in E\otimes_\Zee \Lambda$. Then the
fiber over $x$ of $T_{\mathcal{V}}$ is the divisor $\sum _{i=1}^n\mu_i(e)$. \qed
\end{lemma}

\begin{corollary}\label{morphism} Suppose that $V$ is an irreducible
representation, with lowest weight $\mu$. Let $\hat{T}_0 =  (E\otimes_\Zee
\Lambda)/W_0$ be as in {\rm \S2.1}. Then there is a component
$T_0$ of $T_{\mathcal{V}}$, which has multiplicity one in $T_{\mathcal{V}}$, and
a finite birational morphism
$\varphi \colon \hat{T}_0
\to T_0$ covering the identity on $\Mod$, where the composite map $\hat{T}_0 \to
T_0 \to \Mod \times E_{\rm reg}$ is $(\hat \nu, r)$. Thus
$\hat{T}_0$ is identified with the normalization of a component $T_0$ of
$T_{\mathcal{V}}$.
\end{corollary}
\begin{proof}  The reduction of the spectral cover $T_{\mathcal{V}}\subseteq
\Mod\times E$ is the set of all points $(x,p)$ such that there exists a
lift of $x$ to $e\in E_{\textrm{reg}}\otimes_\Zee\Lambda$ and a weight $\mu_i$
such that 
$\mu_i(e) = p$. Define the map
$\widetilde \varphi\colon E_{\textrm{reg}}\otimes_\Zee \Lambda\to
T_{\mathcal{V}}$ by
$\widetilde \varphi(e) = (\widetilde\nu(e), \mu(e))$, where $\widetilde \nu
\colon E_{\textrm{reg}}\otimes_\Zee\Lambda \to \Mod$ is the quotient morphism.
Clearly,
$\widetilde
\varphi$ descends to a morphism
$\varphi\colon \hat{T}_0\to T_{\mathcal{V}}$, covering the projection to $\Mod$, whose
image is set-theoretically a component
$T_0$ of $T_{\mathcal{V}}$.  Since all weights in $W\cdot \mu$ have multiplicity
one,
$T_0$ is generically reduced and hence is reduced. It will suffice
to show that
$\varphi$ has degree one. Since both
$\hat{T}_0$ and $T_0$ map to $\Mod$ with degree $n$, this is clear.
\end{proof}

\subsection{Some lemmas on weights}

In this section, our goal is to analyze the morphism $\hat{T}_0\to T_0$ of
Corollary~\ref{morphism} more carefully. To do so, we shall have to analyze the
weights in the $W$-orbit of $\mu$. Given two weights $\mu_1, \mu_2$, define the
divisor
$D(\mu_1- \mu_2) \subseteq E_{\rm reg}\otimes_\Zee \Lambda$ via
$$D(\mu_1- \mu_2) = \{e\in E_{\rm reg}\otimes_\Zee \Lambda: \mu_1(e)
=\mu_2(e)\}.$$ More generally, for any weight $\lambda$, let $D(\lambda)=\Ker
\{\lambda\colon E_{\rm reg}\otimes_\Zee \Lambda \to E_{\rm reg}\}$. 

\begin{proposition}\label{minwt} Let $V$ be an irreducible representation of $G$
with lowest weight $\mu$. Suppose that $\mu_1,\mu_2\in W\cdot \mu$ are distinct
weights of
$V$ in the $W$-orbit of $\mu$.  Then exactly one of the following holds:
\begin{enumerate}
\item[\rm (i)] $\mu_1-\mu_2=k\alpha$ for some root $\alpha$ and $k\in \Zee$. In
this case,
$\mu_2$ is the unique element of $W\cdot \mu$ such that $\mu_1-\mu_2$ is a
nonzero real multiple of $\alpha$.
\item[\rm (ii)] There is  an open dense subset of $D(\mu_1- \mu_2)$ on which no
root vanishes.
\end{enumerate}
\end{proposition}
\begin{proof} First suppose that $\mu_1-\mu_2 =k\alpha$ for some root
$\alpha$. Since the difference of two weights of an irreducible
representation is a sum of roots, and the roots are primitive elements in
the root lattice, $k\in \Zee$. The affine line $\mu_1 +t\alpha$ meets the orbit
$W\cdot \mu$ at least at the two points $\mu_1, \mu_2$. Since every element of
$W\cdot\mu$ has the same length (with respect to any Weyl invariant inner
product) and an affine line can meet a sphere in at most two points, we have
proved (i).

Conversely, suppose that $\mu_1-\mu_2$ is not a multiple of any root. Then,
for every root $\alpha$, 
$\Ker (\mu_1-\mu_2)$ and $\Ker \alpha$ are distinct hyperplanes in $\frak h$.
It follows that, for every root $\alpha$, $D(\alpha)$ meets $D(\mu_1- \mu_2)$
in a proper subset. Since there are only finitely many roots, (ii) follows.
\end{proof}

Let us analyze these two possibilities separately. Before we do so, let us
introduce the following notation: 

\begin{defn}\label{notation} For $e\in E_{\textrm{reg}}\otimes_\Zee \Lambda$, let
$\xi_e$ be the  regular $G$-bundle whose S-equivalence class is the image of $e$
in
$\Mod$. Let
$\mathcal{V}_e$ be the vector bundle $\xi_e\times _GV$. 
\end{defn}

\begin{lemma}\label{lemma1} In the above notation, suppose that $\mu_1-\mu_2$
is not a multiple of a root. Then there is an open dense subset $U$ of $D(\mu_1-
\mu_2)$ such that, for all $e\in U$, the following hold:
\begin{enumerate}
\item[\rm (i)] $\operatorname{Stab}_W(e) =\{1\}$.
\item[\rm (ii)] The multiplicity of every point of the divisor $\sum _{w\in
W/W_0}(w\mu)(e)$ is at most two.
\item[\rm (iii)]  If $\xi'$ is the  $H$-bundle defined by $e$, then
$\xi_e=\xi'\times_HG$, and the vector bundle $\mathcal{V}_e$ defined in
Definition~\ref{notation} is a direct sum of line bundles.
\end{enumerate}
\end{lemma}
\begin{proof}  An element $w\in W$ fixes a codimension one subset of
$E_{\textrm{reg}}\otimes_\Zee \Lambda$ if and only if $w$ fixes a codimension one
subset of $\frak h$ if and only if $w=r_\alpha$ is the reflection in a root
$\alpha$. In this case, since $r_\alpha(e) = e-\alpha(e)\alpha\spcheck$ and
$\alpha\spcheck$ is primitive, the fixed set of
$r_\alpha$ is
$D( \alpha)$. Since by hypotheses $D(\alpha) \neq D(\mu_1-\mu_2)$,  (i) 
follows for an open dense subset of
$D(\mu_1- \mu_2)$.

To prove (ii), first suppose that $\mu', \mu''\in W\cdot \mu$ and
that $\mu'-\mu''$ vanishes along a component of $D(\mu_1- \mu_2)$. Then, in
$\frak h$,
$\Ker(\mu'-\mu'') = \Ker (\mu_1- \mu_2)$. Hence $\mu'-\mu''$ is a multiple of
$\mu_1-\mu_2$. Thus, if $\mu', \mu'', \mu'''$ are three distinct
elements of $W\cdot
\mu$ and $\mu ', \mu'', \mu'''$ all agree along a component of $D(\mu_1- \mu_2)$, 
then  $\mu'-\mu''$ and $\mu''-\mu'''$ are both multiples of
$\mu_1-\mu_2$, and so  $\mu', \mu'', \mu'''$ are all contained in an affine line
in $\frak h_\Ar$. But this is impossible since, as $\mu', \mu'', \mu'''$ are
all conjugate via $W$, they all have the same length, and a sphere meets an
affine line in at most two points.

To see (iii), it is clear that $\xi_e$ is regular  since the regular
bundle corresponding to $e$ fails to be split if and only if a root vanishes
on $e$
(by \cite[6.2]{FMI}, in case $E$ is smooth, and by the direct descriptions of
\S1.1 in case $E$ is singular). The  vector bundle
$\mathcal{V}_e$ is thus a direct sum of line bundles.
\end{proof}

The following handles the case where $\mu_1-\mu_2$ is a multiple of a root:

\begin{lemma}\label{lemma2} In the above notation, suppose that
$\mu_1-\mu_2=k\alpha$ for some root $\alpha$ and some $k\in \Zee$, $k\neq 0$.
Then  there is an open dense subset
$U$ of $D(\mu_1- \mu_2)=D(k\alpha)$ such that, for all
$e\in U$, the following hold:
\begin{enumerate}
\item[\rm (i)] If $\mu', \mu''\in W\cdot\mu$ and $\mu'(e) =\mu''(e)$, then
$\mu'-\mu''=\ell\alpha$ for some $\ell \in \Zee$.
\item[\rm (ii)] The multiplicity of every point of the divisor $\sum _{w\in
W/W_0}(w\mu)(e)$ is at most two.
\item[\rm (iii)] $\operatorname{Stab}_W(e) =\{1, r_\alpha\}$ if $e\in
D(\alpha) \subseteq D(k\alpha)$ and $\operatorname{Stab}_W(e) =\{1\}$ if $e\in
D(k\alpha) \setminus D(\alpha)$. Moreover,
$r_\alpha(\mu') = \mu''$ for each pair of distinct weights $\mu', \mu''\in W\cdot
\mu$ such that 
$\mu'(e) =\mu''(e)$.
\item[\rm (iv)] If $e\in
D(k\alpha) \setminus D(\alpha)$, then $\xi_e$ is split, and  $\mathcal{V}_e$ is a
direct sum of line bundles.
\item[\rm (v)]   Suppose that
$e\in D(\alpha)$, and let $G_\alpha$ be the connected
subgroup of
$G$ whose Lie algebra is $\frak g^\alpha \oplus \frak g^{-\alpha}
\oplus \frak h$. Then $G_\alpha$ is isomorphic to $SL_2\times H'$ or to
$SL_2\times _{\Zee/2\Zee}H'$, where $H'$ is the connected subtorus of $H$ whose
Lie algebra is $\Ker\alpha$. If  $\xi_e$ is the corresponding regular
bundle, then $\xi_e =\xi'\times_{(SL_2\times H')}G$, where 
$\xi'=\xi_1\widehat{\boxtimes}\xi_2$ is the bundle over
$SL_2\times H'$ such that $\xi_1$ is the principal $SL_2$-bundle
corresponding to $I_2$ and $\xi_2$ is the $H'$-bundle such that
$\xi'\times _{H'}G$ is   the split bundle whose image in
$\Mod$ is
$e$. 
\end{enumerate}
\end{lemma}
\begin{proof}  If $\mu', \mu''\in
W\cdot\mu$ and
$\mu'-\mu''$ vanishes on a nonempty open subset of $D(\mu_1- \mu_2)$, then
$\mu'-\mu'' = t(\mu_1-\mu_2) =\ell\alpha$ for some real number $\ell$. By the
proof of Proposition~\ref{minwt}, $\ell \in \Zee$, proving (i). Part (ii) and 
the first part of (iii) follow  as in the proof of the previous lemma. The second
part follows from the fact
\cite[p.\ 124, Prop.\ 3(i)]{Bour78} that
$\mu'-j\alpha$ is a weight whose length is the same as that of $\mu'$ if and only
if $j= \mu'(\alpha\spcheck)$, and in this case by definition
$\mu'-j\alpha= r_\alpha(\mu')$. Parts (iv) and (v) follow from the explicit
description of the regular representative $\xi_e$ given in \cite{FMI} in the
smooth case, and by the discussion of \S1.1 in the singular case.
\end{proof}

We can now describe the morphism $\varphi\colon \hat{T}_0 \to T_0\subseteq
T_{\mathcal{V}}$ in more detail. Recall that $\varphi$ is the morphism induced
from $\widetilde \varphi\colon E_{\textrm{reg}}\otimes _\Zee\Lambda \to
T_0\subseteq \Mod \times E$ defined by $\widetilde \varphi(e) = (\widetilde\nu(e),
\mu(e))$, where $\widetilde \nu\colon E_{\textrm{reg}}\otimes _\Zee\Lambda \to
\Mod$ is the quotient map.

\begin{lemma}\label{singofmap} In the above notation, suppose that $e\in
E_{\rm reg}\otimes_\Zee \Lambda$ , and let $y\in \hat{T}_0$ be the image of $e$.
\begin{enumerate}
\item[\rm (i)] If, for all   $\mu'\in W\cdot \mu$ such that $\mu'\neq \mu$, 
$e\notin D(\mu-\mu')$, then $\varphi^{-1}(\varphi(y))
=\{y\}$, and the morphism
$\varphi\colon \hat{T}_0\to T_0$ is an isomorphism in a neighborhood of $y$.
\item[\rm (ii)] If   $\mu'\in W\cdot \mu$, $\mu'\neq \mu$, and $e$ is a
generic point of $D(\mu-\mu')$, where either
$\mu-\mu'$ is not a multiple of a root, or
$\mu-\mu'=k\alpha$ for some root $\alpha$ and $e\in D(k\alpha) \setminus
D(\alpha)$, then there is an open neighborhood $U$ in $T_0$ of $\varphi(y)$ such
that
$\varphi^{-1}(U) = U_1\cup U_2$, $\varphi|U_i$ is an isomorphism onto a smooth
divisor in $\Mod \times E$, and $\varphi(U_1)$ meets $\varphi(U_2)$
transversally. 
\item[\rm (iii)] If   $\mu'\in W\cdot \mu$, $\mu'\neq \mu$, and $e$ is a generic
point of $D(\mu-\mu')$, where $\mu-\mu'=k\alpha$ for some root $\alpha$ and
$e\in D(\alpha)$, then $\varphi\colon
\hat{T}_0\to T_0$ is an isomorphism in a neighborhood of $y$.
\end{enumerate}
\end{lemma}
\begin{proof} Suppose that $e\notin D(\mu-\mu')$ for every $\mu'\in W\cdot \mu$
such that
$\mu'\neq \mu$. Given $y'\in
\hat{T}_0$ which is the image of  $e'\in E_{\textrm{reg}}\otimes_\Zee \Lambda$,
suppose that
$\varphi(y) =\varphi(y')$. Then $e'=w(e)$ for some $w\in W$ and $\mu(e') =
\mu(we) =\mu(e)$. It follows that $w^*\mu = \mu$, so that  $w\in
W_0$, i.e.\ 
$e$ and $e'$ have the same image $y\in \hat{T}_0$. Thus $\varphi^{-1}(\varphi(y))
=\{y\}$.  If $w\in \operatorname{Stab}_W(e)$, then  
$(w\mu)(e)
= \mu(e)$, and so $w\mu =\mu$, i.e.\ $w \in W_0$. Thus
$\operatorname{Stab}_W(e) = \operatorname{Stab}_{W_0}(e)$, so that the
morphism $\hat\nu \colon \hat{T}_0 \to \Mod$ is a local diffeomorphism near $y$.
Since
$\hat\nu =\pi_1\circ \varphi$, the morphism $\varphi$ must also be a local
diffeomorphism onto its image near $y$. This proves (i).

Now suppose that the hypotheses of (ii) hold and let $e'$ and $y'$ be as above. By
(ii) of Lemma~\ref{lemma1} or (ii) of Lemma~\ref{lemma2}, if $e$ is generic, then
$\mu'$ is the unique element of $W\cdot \mu$ such that $\mu'(e) = \mu(e)$.
Suppose that
$\mu'=w\mu$. If
$\varphi(y')=\varphi(y)$, then necessarily $e' = w^{-1}(e)$, and so
$\varphi^{-1}(\varphi(y)) =\{y,y'\}$. Moreover, since $\operatorname{Stab}_W(e)
=\{1\}$, $y\neq y'$ and the map $\hat{T}_0 \to \Mod$ is a local diffeomorphism onto
its image at both
$y$ and $y'$. Hence the same is true for $\varphi$. Since $w\mu =\mu'\neq \mu$, it
is easy to see that the two divisors meet transversally at $\varphi(y)$. This
proves (ii).

Finally suppose that we are in Case (iii). In this case $\operatorname{Stab}_W(e)
=\{1, r_\alpha\}$ and $\mu' =r_\alpha(\mu)$. Thus $\operatorname{Stab}_{W_0}(e)
=\{1\}$. It follows as in the previous case
that, if $\varphi(y')=\varphi(y)$, then $e' = r_\alpha(e) = e$. Thus
$\varphi^{-1}(\varphi(y)) =\{y\}$. 
For $e$ generic, we can identify the tangent space to $\hat{T}_0$ at $y$  with
$\frak h$ and the kernel of the differential of the map from $\hat{T}_0$ to
$\Mod$ can be identified with $\Cee\cdot \alpha\spcheck \subseteq \frak h$. Since
$\mu (\alpha\spcheck)
\neq 0$, it follows that the differential of $\varphi$ at $y$ is injective. This
proves (iii).
\end{proof}

\subsection{The case of a minuscule or quasiminuscule representation} 

We now apply the above results in case the representation $\rho$ is minuscule or
quasiminuscule. Given $\rho$, for each pair of integers $(a,b)$, we have
constructed a vector bundle $\mathcal{V}_{a,b}$ in  \S2.2, and the
spectral cover $T_{\mathcal{V}_{a,b}}$ is independent of $(a,b)$. As before, we
fix some pair $(a,b)$ and set $\mathcal{V}=\mathcal{V}_{a,b}$.

\begin{theorem}\label{spectralcoverthm} Let $E$ be singular, and suppose
that $\rho$ is minuscule. 
Then there exists a line bundle $M$ on
$\hat{T}_0$ such that
$\mathcal{V} = (\hat\nu\times \Id)_*(\pi_1^*M \otimes (r\times \Id)^*\Poin)$.
A similar result holds if $E$ is smooth, provided that we replace $\hat\nu\colon
\hat{T}_0\to \Mod$ by $\hat\nu\colon (\hat{T}_0)_{\rm reg} \to \Mod_{\rm reg}$. 
\end{theorem}

\begin{proof} In the minuscule case, $T_0=T_{\mathcal{V}}$.
Note that the pullback of the Poincar\'e line bundle on $T_{\mathcal{V}}\times E$
to
$\hat{T}_0$ is simply
$(r\times \Id)^*\Poin$. Thus, it suffices to show that the normalization morphism
$\varphi\colon \hat{T}_0
\to T_{\mathcal{V}}$ and the bundle $\mathcal{V}$ satisfy the hypotheses of
Proposition~\ref{normalization}. Let $D_0$   be the union in $\Mod$ of the
images of the divisors $D(\mu-\mu')$, and let $D_0'$ be the dense open subset
of $D_0$ implicitly defined by Lemmas~\ref{lemma1} and  \ref{lemma2}. Note that,
since
$\rho$ is minuscule,  if the difference $\mu-\mu'$ of two weights is a
multiple of a root $\alpha$, then in fact this multiple is $\pm 1$, so that Case
(iv) of Lemma~\ref{lemma2} does not arise.

If $x\in \Mod$ does not lie in $D_0$, then $T_{\mathcal{V}}$ is smooth over $x$
and
$\mathcal{V}_x$ is regular. If $x\in D_0'$ lies under a point
of $D(\mu-\mu')$ satisfying the hypotheses of Lemma~\ref{lemma1}, then
$\mathcal{V}_x$ satisfies (ii)(a) (with $\mathcal{V}_x(q)$ of rank one) or
(ii)(b) of Proposition~\ref{normalization}. Otherwise, we are in case (v) of
Lemma~\ref{lemma2}, and $x\in \Mod$ lies under a generic point $e\in D(\alpha)$
for some root $\alpha$.  In this case, there is a connected subgroup
$G_\alpha$ of
$G$ isomorphic to $SL_2\times H'$ or to $SL_2\times _{\Zee/2\Zee}H'$, the
structure group of $\xi_e$ reduces to
$G_\alpha$ and lifts to $SL_2\times H'$, and on the $SL_2$ factor it is the
principal bundle corresponding to 
$I_2$. Under the induced action of $SL_2$ on $V$, $V$ decomposes into 
$\bigoplus _iW_i\oplus \bigoplus _jU_j$, where the $W_i$ are irreducible,
nontrivial representations of
$SL_2$, and the $U_i$ are the trivial representation. Since $\rho$ is minuscule,
$\dim W_i=2$ by e.g. \cite[p.\ 128, Prop.\ 7]{Bour78}. Thus, for the action of
$SL_2\times H'$ on $V$, $V$ decomposes into 
$\bigoplus _i(W_i\otimes \chi_i)\oplus \bigoplus _j(U_j\otimes \chi_j)$, where
the
$\chi_i$ and $\chi_j$ are characters of $H'$. It then follows that the vector
bundle $\mathcal{V}_e$ associated to
$\xi_e$ is of the form $\bigoplus _i(I_2\otimes
\lambda_i)\oplus \bigoplus _j\lambda_j$. Here the
$\lambda_i$ in the first summand are the line bundles of degree zero
corresponding to
$\mu'(e)$ for
$\mu'\in W\cdot \mu$ such that $\mu' -\mu''=\pm\alpha$ for some $\mu''\in
W\cdot \mu$, and the $\lambda_j$ in the second summand are the line bundles of 
degree zero corresponding to $\mu'(e)$ for
$\mu'\in W\cdot \mu$ such that $\mu'(e)\neq \mu''(e)$ if $\mu''\in W\cdot
\mu$, $\mu''\neq \mu'$. In particular, the line bundles corresponding to
the various summands are all distinct. Thus, the vector bundle associated to
$\xi_e$ is regular.
In this case,
$\varphi$ fulfills the hypotheses of (ii)(a) of Proposition~\ref{normalization},
completing the proof of the theorem.
\end{proof}

\begin{remark} It is an interesting problem to identify the line bundle $M$ of
Theorem~\ref{spectralcoverthm} in case $E$ is smooth, where it will depend on the
choice of $(a,b)$. For groups of type
$A_n$ or $C_n$, the answer is essentially contained in \cite[Corollary 3.4]{FMW}.
\end{remark} 

Now we turn to quasiminuscule representations.

\begin{theorem}\label{spectralcoverthm2} Let $E$ be singular, and suppose
that $\rho$ is quasiminuscule but not minuscule. 
Define the sheaf $\overline{\mathcal{V}}$ by the exact sequence
$$0 \to \pi_1^*\pi_1{}_*\mathcal{V} \to \mathcal{V} \to \overline{\mathcal{V}} 
\to 0.$$ Then $\overline{\mathcal{V}}$ is a vector bundle, and its associated
spectral cover is
$T_0$, the image of $\hat{T}_0$. Moreover, there exists a line bundle
$M$ on
$\hat{T}_0$ such that
$$\overline{\mathcal{V}} = (\hat\nu\times \Id)_*(\pi_1^*M \otimes (r\times
\Id)^*\Poin).$$ A similar result holds if $E$ is smooth, provided that we replace
$\hat\nu\colon
\hat{T}_0\to \Mod$ by $\hat\nu\colon (\hat{T}_0)_{\rm reg} \to \Mod_{\rm reg}$. 
\end{theorem}
\begin{proof} Suppose that the multiplicity of the trivial weight in $\rho$ is
$s$. Let $e\in E_{\textrm{reg}}\otimes_\Zee \Lambda$ and let $\xi_e$ be the
regular vector bundle whose S-equivalence class corresponds  to the image of $e$
in $\Mod$. We first claim that $h^0(E; \xi_e\times _GV) = s$ for every $e$. Of
course, if $G$ is simply laced, then $\rho$ is the adjoint representation, $s=r$,
and the claim is essentially the definition of a regular bundle. 
In case $G$ is not simply laced, the nonzero weights are the short roots. Let
$e\in E_{\textrm{reg}}\otimes_\Zee \Lambda$ and let $R_e$ be the set of roots
which are trivial on
$e$. Let $X$ be a principal nilpotent element for the reductive Lie algebra
$\frak z_{\frak g}(e)$ which is generated by $\frak h$ and the root spaces $\frak
g^\alpha$, $\alpha \in R_e$.  By the recipe for constructing
$\xi_e$ given in \cite{FMI} in case $E$ is smooth, and by the discussion in \S1.1
in case $E$ is singular,
$H^0(E;\xi_e\times _GV)$ is given by the kernel of
$\rho_*(X)$ acting on $V_e$, where $V_e$ is the sum of the weight spaces
(including those for the trivial weight) which annihilate $e$. The trivial weight
contributes a subspace of $V_e$ of dimension
$s$ and the remaining weight spaces annihilating $e$ are one-dimensional
subspaces of $V_e$  corresponding to the short roots in $R_e$. We may assume that
we have chosen $X$ so that there is an
$h_e\in \frak h \subseteq \frak z_{\frak g}(e)$ such that $(h_e, X)$ completes to
an $\frak{sl}_2$-triple in $\frak z_{\frak g}(e)$. Thus $(\rho_*(h_e), \rho_*(X))$
completes to some
$\frak{sl}_2$-triple in $\operatorname{End} V$. Since
$X$ is principal in
$\frak z_{\frak g}(e)$,  for all $\alpha\in R_e$, $\alpha(h_e)$ is an even
nonzero integer. Thus the
kernel of $\rho_*(h_e)$ is the weight zero subspace of $V$ and so $\dim
\Ker\{\rho_*(h_e) \colon V_e\to V_e\} = s$. Since all of the eigenvalues of
$\rho_*(h_e)$ are even, it follows from the  classification of
representations of
$\frak{sl}_2$  that
$\dim \Ker\{\rho_*(X)
\colon V_e\to V_e\} = s$ as well.

Since $h^0(E; \xi_e\times _GV)$ is constant, the direct image
$\pi_1{}_*\mathcal{V}$ is locally free of rank $s$ and the induced map
$\pi_1{}_*\mathcal{V}/\frak m_x\pi_1{}_*\mathcal{V} \to H^0(E; \mathcal{V}_x)$ is
an isomorphism for every $x\in \Mod$. By semistability, the map $H^0(E;
\mathcal{V}_x) \otimes \scrO_E\to \mathcal{V}_x$ is injective for every $x\in
\Mod$. It follows that the induced homomorphism $\pi_1^*\pi_1{}_*\mathcal{V} \to
\mathcal{V}$ is injective on every fiber, and so the cokernel 
$\overline{\mathcal{V}}$ is a subbundle as claimed. Clearly the spectral cover
corresponding to
$\pi_1^*\pi_1{}_*\mathcal{V}$ is just the trivial section of $\Mod\times E$ with
multiplicity $s$, and hence $T_{\overline{\mathcal{V}}} = T_0$.

As in the proof of Theorem~\ref{spectralcoverthm}, we take $D_0$ to be the union
in
$\Mod$ of the images of the divisors $D(\mu-\mu')$, and $D_0'$ to be the dense
open subset of
$D_0$ implicitly defined by Lemmas~\ref{lemma1} and  \ref{lemma2}. To complete the
proof, we must analyze the corresponding bundles
$\overline{\mathcal{V}}_x$, where either $x\notin D_0$ or $x\in D_0'$. The
arguments are very similar to those given in the proof of
Theorem~\ref{spectralcoverthm}, and we shall be brief. The only essentially new
case is when a weight, i.e.\ a short root $\alpha$, vanishes. In this case, for a
generic $x\in D(\alpha)$, there is a unique factor of
$\mathcal{V}_x$ isomorphic to $I_3$, and the remaining summands of
$\mathcal{V}_x$ are of the form $I_2\otimes \lambda_i$ or $\lambda_j$, where the
$\lambda_i, \lambda_j$ are pairwise distinct nontrivial line bundles, together
with $s-1$ copies of $\scrO_E$. In this case, $\overline{\mathcal{V}}_x$ has
exactly the same  summands, except that the $I_3$ factor is replaced by an $I_2$
and the remaining $s-1$ copies of $\scrO_E$ are no longer present. Hence
$\overline{\mathcal{V}}_x$ is regular.
\end{proof}

\subsection{Some remarks on representations}

Given a regular $G$-bundle $\xi$ and an irreducible representation $\rho\colon
G\to GL(V)$ of
$G$, when is the associated vector bundle $\xi\times _GV$ regular? In this
section, we state some results along these lines. The proofs are straightforward
and will not be given.

\begin{proposition}\label{prop1} Let $\rho\colon
G\to GL(V)$ be an irreducible representation. Then the following are equivalent:
\begin{enumerate}
\item[\rm (i)] For every smooth elliptic curve $E$ and every regular $G$-bundle
$\xi$ over $E$, the associated vector bundle $\xi\times _GV$ is a regular vector
bundle.
\item[\rm (ii)] For every regular element $g\in G$, $\rho(g)$ is a regular
element of $GL(V)$.
\item[\rm (iii)] The multiplicity of every weight of $\rho$ is one, and for
every pair $\mu_1,\mu_2$ of distinct weights of $\rho$, $\mu_1-\mu_2$ is a root.
\item[\rm (iv)] Either $G=SL_n$ and $\rho$ is the standard representation or its
dual, or $G=Sp(2n)$ and $\rho$ is the standard representation. \qed
\end{enumerate}
\end{proposition}

A weaker but related question is the following: let $\xi$ be the regular
$G$-bundle S-equivalent to the trivial $G$-bundle. When is $\xi\times _GV$
regular? The answer is as follows:

\begin{proposition}\label{prop2} Let $\rho\colon
G\to GL(V)$ be an irreducible representation. Then the following are equivalent:
\begin{enumerate}
\item[\rm (i)] For every smooth elliptic curve $E$ and every regular $G$-bundle
$\xi$ over $E$ S-equivalent to the trivial bundle, the associated vector bundle
$\xi\times _GV$ is a regular vector bundle.
\item[\rm (ii)] For every regular element $X\in \frak g$, $\rho_*(X)$ is a
regular element of $\operatorname{End} V$.
\item[\rm (iii)] If $X$ is a principal nilpotent element in $\frak g$, then
$\rho_*(X)$ is a regular element of $\operatorname{End} V$.
\item[\rm (iv)] The multiplicity of every weight of $\rho$ is one, and for
every pair $\mu_1,\mu_2$ of   weights of $\rho$, $\mu_1-\mu_2$ is a multiple
of a root.
\item[\rm (v)] The group $G$ and the representation $\rho$ are as in {\rm (iv)}
of Proposition {\rm\ref{prop1}}, or $G=SL_2$ and $\rho$ is any irreducible
representation, or $G$ is of type $B_n$ or $G_2$ and $\rho$ is the unique
quasiminuscule, non-minuscule representation of $G$.
\qed
\end{enumerate}
\end{proposition}

\section{Kostant and Steinberg sections and the proof of the main theorem}

\subsection{Conjugation into the slice}

Every $\ad G$-invariant function on $\frak g$ may be identified with a
$W$-invariant function on $\frak h$. In this way, the adjoint quotient
of $\frak 
g$ is identified with $\frak h/W$. Let $p\colon \frak g\to \frak h/W$ be the
induced morphism, and let $\frak g_{\textrm{reg}}$ be the open dense subset of
$\frak g$ consisting of regular elements. In \cite{Kostant}, Kostant
constructed 
a section of $p$ whose image lies in $\frak g_{\textrm{reg}}$ as follows. Let
$\frak b$ be a Borel subalgebra containing $\frak h$, with nilpotent radical
$\frak u$. Thus, there is a unique  set of positive roots $R^+$ for
$(\frak g, \frak h)$ such that $\frak u =\bigoplus_{\alpha \in R^+}\frak
g^\alpha$. Fix a principal $\frak{sl}_2$-triple $(X, h_0, Y)$ adapted to $\frak
b$. By definition, this means that $X\in \frak u$ is a principal nilpotent
element of $\frak g$, $[h_0, X] =2X$, $[h_0, Y] =-2Y$, and $[X,Y] =h_0\in \frak
h$. It follows automatically that $Y\in \frak u_- =  \bigoplus_{-\alpha \in
R^+}\frak g^\alpha$.   Kostant  has proved the following \cite[Theorem
0.10]{Kostant}:

\begin{theorem} Let $Z$ be a vector subspace of $\frak g$ which is
complementary  
to $\operatorname{Im}(\ad X)$ and is invariant under $\ad h_0$. Then the affine
space $Z+X$ is contained in $\frak g_{\rm reg}$ and the morphism $Z+X \to
\frak h/W$ induced by the projection $\frak g \to \frak h/W$ is an isomorphism.
\qed
\end{theorem}

We call the inverse morphism $\ov\Sigma_Z\colon \frak h/W \to Z+X \subseteq \frak
g_{\textrm{reg}}$ the 
\textsl{Kostant section} defined by $Z$. 
Denote the 
composition $\frak h\to \frak h/W \to Z+X \to \frak g_{\textrm{reg}} \subseteq
\frak g$ by $\Sigma_Z \colon \frak h \to \frak g$. Clearly,
$\Sigma_Z$ is a $W$-invariant morphism from $\frak h$ to $\frak
g$. 

We would like to compare the $\ov\Sigma_Z$ for different choices of $Z$,
and would also like to compare $\Sigma_Z$
with the inclusion $i\colon \frak
h\to \frak 
g$. Of course, $\Sigma_Z$  can never be 
conjugate to the inclusion, not even pointwise, 
because not all elements of $\frak h$ are 
regular in $\frak g$. According to the next result,  for every choice of
$Z$, if we replace the inclusion morphism $i$ by
$i+X$, then the resulting morphism is indeed conjugate to $\Sigma_Z$   via a
morphism $\frak h\to G$:

\begin{theorem}\label{conj} Let  $(X, h_0, Y)$ be a principal
$\frak{sl}_2$-triple 
adapted to
$\frak b$ and let $Z$ and $Z'$ be linear complements to $\operatorname{Im}(\ad X)$
 invariant under $\ad h_0$. Let $\ov\Sigma_Z$ and $\ov\Sigma_{Z'}$ be the Kostant
sections defined by $Z$ and $Z'$ respectively. Let
$i+X\colon
\frak h\to \frak g$ be the  map $h\mapsto h+X$. Then:
\begin{enumerate}
\item[\rm (i)] $\ov\Sigma_Z$ and $\ov\Sigma_{Z'}$ are conjugate under the
adjoint action of $G$ on $\frak g$, i.e.\ there exists a morphism
$\psi \colon \frak h/W\to G$ such that
$$\Ad (\psi)\circ \ov\Sigma_Z=\ov\Sigma_{Z'}.$$ 
\item[\rm (ii)] The map $i+X$ is an embedding of $\frak h$ into $\frak
g_{\textrm{reg}}$.
\item[\rm (iii)] For each $h\in \frak h$, $p(h+X) = p(h)$, so that
$h+X$ and $h$ 
are identified in the adjoint quotient. More precisely, there is a
one-parameter 
subgroup $\beta\colon \Cee^*\to G$ such that, for all $h\in \frak h$,
$$\lim_{t\to 0}\Ad \beta(t)(h+X) = h.$$
\item[\rm (iv)] There exists a  morphism $\Psi\colon \frak h\to \frak u_-$
such that, for all $h\in H$,
$$\Ad (\exp(\Psi(h)))(h+X) = \Sigma_Z (h).$$
\item[\rm (v)] The morphism $\Psi$ is the unique morphism from $\frak h$ to
$\frak u_-$ with the property that, for all
$h\in \frak h$, 
$$\Ad (\exp(\Psi(h)))(h+X) \in Z+X.$$ 
\end{enumerate}
\end{theorem}

By (i), all Kostant sections are equivalent under the adjoint
action of $G$ on $\frak g$, a result which is subsumed in the main result of
\cite{Auto}. By (iv), after making the base change $\frak h\to \frak
h/W$, we can  
conjugate any Kostant section until its image lies in a Borel
subalgebra. This is 
of course not possible for the  Kostant section itself, since
otherwise we would be able to find a section to the morphism $\frak h\to
\frak h/W$.  

We will deduce Theorem~\ref{conj} from a more general result,  
proved in \S\ref{nextsection}.

\subsection{More general $\frak{sl}_2$-triples}\label{nextsection}

For the moment, let $(X,h_0,Y)$ denote an arbitrary $\frak{sl}_2$-triple, not
necessarily principal or adapted to $\frak b$. The endomorphism $\ad h_0$ of
$\frak g$ is semisimple and its eigenvalues are integers. Let $\frak g_{\ell}$ be
the $\ell$-eigenspace of $\ad h_0$. Thus $\frak g$ is graded and the bracket
is compatible with the grading. Moreover, $X\in \frak g_2$. Let $\frak
g_-$ be the nilpotent subalgebra $\bigoplus_{\ell < 0}\frak g_{\ell}$, let $G_-$
be the corresponding connected unipotent subgroup of $G$,  and let
$\frak g_{\leq 0} =\bigoplus_{\ell \leq 0}\frak g_{\ell}$. 
By the \textsl{degree} $d(x)$ of an element $x\in \frak g$, we mean the largest
value of
$\ell$ such that the component of $x$ in $\frak g_\ell$ is nonzero. By the
classification of finite-dimensional representations of $\frak{sl}_2$, $\ad
X\colon \frak g_\ell \to \frak g_{\ell +2}$ is injective if $\ell < 0$ and
surjective if $\ell \geq -1$.

\begin{proposition}\label{main} Let $(X,h_0,Y)$ be an $\frak{sl}_2$-triple. Let
$Z\subseteq
\frak g_{\leq 0}$ be a complementary linear subspace to the image of $\ad X$
invariant under $h_0$. Let $\tau \in \frak g_{\leq 0}$. Then there is a unique
$\psi \in \frak g_-$ such that $\Ad (\exp(\psi))(\tau + X) \in Z+X$.
\end{proposition}
\begin{proof} Write $\tau = \tau_0+\ad X(\tau_1)$, where $\tau_0\in Z$
and we may 
assume that $\tau_1\in \frak g_-$, in fact that $d(\tau_1)\leq
-2$. The proof of existence of $\psi$ is by induction on
$k=d(\tau_1)\leq -2$. If $k$ is 
sufficiently negative, then $\tau_1=0$ and we can take
$\psi=0$. Suppose that we 
know the result for all $\tau$ such that $d(\tau_1) \leq k-1$, and let
$\tau$ be 
an element such that $d(\tau_1) = k$. Let $\psi' = \tau_1$. Then
\begin{align*}
\Ad(\exp(\psi'))(\tau+X) &=   \exp(\ad \psi')(\tau + X) \\
&= \tau + X
-\ad X(\tau_1) + \sum _{n=1}^\infty\frac{1}{n!}\ad (\tau_1)^n(\tau) + \sum
_{n=2}^\infty\frac{1}{n!}\ad (\tau_1)^n(X).
\end{align*}
We can write this as $\tau_0 + X +\alpha$, where $d(\alpha) \leq d(\tau_1)$.
Let $\alpha = \tau_0' + \ad X(\tau_2)$, where $\tau_0'\in Z$ and
$\tau_2 \in \frak 
g_-$. Since $Z$ is $h_0$-invariant, $d(\tau') \leq d(\alpha)$ and thus
$d(\ad X(\tau_2)) \leq d(\alpha)$. We can thus assume that $d(\tau_2)
\leq d(\alpha) -2
\leq k-2$. Thus 
by induction there exists a $\psi_1$ such that
$\Ad(\exp(\psi_1))(\tau_0+\alpha + X) \in Z+X$. Since $\exp (\psi_1)$ and
$\exp(\psi')$ lie in the unipotent group $G_-$, $\psi =\log (\exp
(\psi_1)\exp(\psi'))$ is defined and clearly $\Ad(\exp(\psi))(\tau + X)
\in Z+X$. This completes the inductive step and proves the existence of $\psi$.

To see the uniqueness, suppose that $\Ad (\exp(\psi_i))(\tau + X) \in Z+X$ for
$i=1,2$. Writing $\log(\exp(\psi_1)\exp(-\psi_2)) = \psi$, we see that
it suffices 
to prove: if there exists $\zeta \in Z$ and $\psi\in \frak g_-$ such that
$\Ad(\exp(\psi))(\zeta + X) \in Z+X$, then $\psi =0$. If $\psi \neq 0$, write
$\psi = \psi_k+ \psi'$, where $\psi_k$ is homogeneous of degree $k<0$ and
$d(\psi')< k$. As before, we write
$$\Ad (\exp(\psi))(\zeta + X) =\zeta + X + [\psi_k, X] + \gamma ,$$
where $d(\gamma) < d([\psi_k, X]) = k+2$. Since $\zeta \in Z$ and $Z$
is a vector 
subspace of $\frak g$, it follows that $[\psi_k, X] + \gamma\in
Z$. Since $Z$ is 
invariant under the action of $h_0$, each homogeneous component of
$[\psi_k, X] + 
\gamma$ lies in $Z$. But as  $d(\gamma) < d([\psi_k, X])$, it follows that
$[\psi_k, X]\in Z$. Since $Z$ is a complement to the image of $\ad X$,
$[\psi_k, 
X] =0$, and hence $\psi_k=0$ since $\psi_k\in \frak g_k$ and
$k<0$. This is a contradiction. Hence $\psi =0$. 
\end{proof}

\begin{remark} Suppose that $\frak g_{-1}=0$, which is the case if $(X,h_0, Y)$
is principal. Then, for all $g\in G_-$, $\Ad(g)(X)-X\in \frak g_{\leq 0}$. It
follows that, for all $g\in G_-$ and $z\in \frak g_-$, $\Ad (g)(z+X)-X \in
\frak g_{\leq 0}$. Thus  
$$(g,z) \mapsto \Ad (g)(z+X)-X$$
defines a morphism
$$G_-\times (Z+X) \to \frak g_{\leq 0}.$$
It follows from Proposition~\ref{main} that this morphism defines an isomorphism
from $G_-\times (Z+X)$ to $\frak g_{\leq 0}$. In this essentially equivalent
form, Proposition~\ref{main} is due to Kostant \cite[Theorem 1.2]{Kostant2}.
\end{remark}

\subsection{A parametrized version}

We continue to assume that $(X,h_0,Y)$ is an arbitrary
 $\frak{sl}_2$-triple. Let 
 $R$ be a finitely generated $\Cee$-algebra. The exponential map $\exp\colon
\frak g_- \to G$ extends to a map $\frak g_-\otimes _\Cee R \to G(R)$, where
$G(R)$ denotes the group of $R$-valued points of $G$, in other words the set of
morphisms $\Spec R \to G$. Viewing $\frak g_-\otimes _\Cee R$ as the set of
morphisms $\Spec R \to \frak g_-$, this map is the composition $\psi\mapsto
\exp\circ \psi\colon \Spec R \to G$. The group $G(R)$ acts on $\frak g\otimes
_\Cee R$ via $\Ad$. We still have the formula $\Ad (\exp(\psi)) = \exp
 (\ad\psi)$. 
From this, it is clear that the proof of Proposition~\ref{main} also proves:

\begin{proposition}\label{param} Let $(X,h_0,Y)$ be an
 $\frak{sl}_2$-triple. Let $Z\subseteq
\frak g_{\leq 0}$ be a complementary linear subspace to the image of $\ad X$
invariant under $h_0$. Let $\tau \in \frak g_{\leq 0}\otimes _\Cee R$. Then there
is a unique
$\psi \in \frak g_-\otimes _\Cee R$ such that $\Ad (\exp(\psi))(\tau + X) \in
Z\otimes _\Cee R +X\subseteq \frak g\otimes _\Cee R$. \qed
\end{proposition}

\begin{remark}\label{homogeneous} An examination of the inductive argument shows
the following: suppose that $R =\Cee[s_1, \dots, s_n]$ is the coordinate ring of
affine $n$-space with the usual grading,   and that   $\tau\in \frak
g_0\otimes_\Cee R_1$, so that the components of $\tau$ are linear polynomials.
Then $\psi =
\sum _{k=1}^N\psi _{-k}$, where $\psi_{-k}\in
\frak g_{-2k}\otimes_\Cee R_{k+1}$ and so the components of $\psi_{-k}$ are
homogeneous polynomials of degree
$k+1$.
\end{remark}

We will also need the following lemma to prove Part (i) of Theorem~\ref{conj}:

\begin{lemma}\label{iso} Let $Z,Z'\subseteq \frak g_-$ be two
linear subspaces, both complementary  to
$\operatorname{Im}(\ad X)$ and  invariant under $h_0$. Then there
exists a morphism
$\gamma \colon Z \to \frak g_-$ with $\gamma(0) =0$ such that the morphism
$$z\mapsto \Ad (\exp(\gamma(z)))(z+X) - X$$
is an isomorphism from $Z$ to $Z'$.
\end{lemma}
\begin{proof} Let $R$ be the coordinate ring of $Z$. Applying
Proposition~\ref{param} to the inclusion of $Z$ in $\frak g_-$, which
corresponds to an element of
$R\otimes_\Cee \frak g_-$, we see that there is a unique morphism
$\gamma\colon Z \to 
\frak g_-$ with $\Ad (\exp(\gamma(z)))(z+X)\in Z'+X$ for all $z\in
Z$. Clearly $\gamma(0) =0$ by uniqueness. Define the 
morphism $\alpha\colon Z\to Z'$ by
$\alpha(z)=
\Ad (\exp(\gamma(z)))(z+X) - X$. Symmetrically, there is a unique
morphism $\gamma'\colon Z' \to \frak g_-$ such that, for all $z'\in Z'$, $\Ad
(\exp(\gamma(z')))(z'+X)\in Z+X$. If we set
$\beta (z') =  \Ad (\exp(\gamma(z')))(z'+X) - X$, then $\beta$ is a
morphism from 
$Z'$ to $Z$. By the uniqueness statements, it follows that
$$\Ad(\exp(\gamma'(\alpha(z))))\circ\Ad( \exp (\gamma(z)))(z+X) = z+X$$
for all $z\in Z$, and hence that $\beta\circ \alpha =\Id$. Symmetrically,
$\alpha\circ \beta = \Id$ as well, so that
$\alpha$ and $\beta$ are inverse to one another.
\end{proof}

\subsection{Proof of Theorem~\ref{conj}}

We   now prove Theorem~\ref{conj}. Let $(X,h_0, Y)$ be a principal
$\frak{sl}_2$-triple adapted to $\frak b$. Since $h_0$ is regular, $\frak g_0 =
\frak h$ and $\frak g_- =\frak u_-$. 
Let $Z$ and $Z'$ be complementary linear subspaces to $\operatorname{Im}(\ad X)$
and invariant under $h_0$. Both $Z$ and $Z'$ lie in $\frak g_-$. By
Lemma~\ref{iso} there is a morphism $\gamma\colon
Z\to \frak g_-$ such that $z+X\mapsto\Ad(\exp(\gamma(z)))(z+X)$ is an
isomorphism  
$I\colon Z+X$ to
$Z'+X$. 
Since $I$ is given by the adjoint action of $G$ on $\frak
g$,
it follows that $I$ is compatible with the projections to the adjoint
quotient $\frak h/W$. In particular, $\Ad(\exp(\gamma))\circ
\ov\Sigma_Z=\ov\Sigma_{Z'}$,  proving Part (i) of Theorem~\ref{conj}.

The map $\Cee\to \frak h$ defined by $a\mapsto ah_0$ exponentiates to a
one-parameter subgroup $\beta\colon \Cee^* \to H$ such that $\Ad\beta(t)(X) =
t^2X$. Thus, for all $h\in H$, $\Ad\beta(t)(h+X) = h+t^2X$, and hence
$\lim_{t\to 
0}\Ad \beta(t)(h+X) = h$. This proves Part (iii) of Theorem~\ref{conj}.

Viewing the inclusion $i\colon \frak h \to
\frak g$ as a morphism from $\frak h$ to $\frak g_0$, it follows from
Proposition~\ref{param} that there is a unique morphism $\Psi\colon
\frak h\to 
\frak g_-$ such that 
$$\Ad(\Psi(h))(h+X) \in Z+X$$ for all $h\in \frak
h$. This 
proves Part (v) of Theorem~\ref{conj}. Note that, by
Remark~\ref{homogeneous}, $\Psi =
 \sum _{k=1}^N\Psi
_{-k}$, where $\Psi_{-k}\in \frak g_{-2k}\otimes_\Cee R$ is a homogeneous
polynomial of degree $k+1$.

Now define $\Sigma'_Z\colon \frak h \to Z+X$ by
$$\Sigma'_Z(h) = \Ad\exp(\Psi(h))(h+X).$$

\begin{lemma} The morphism $\Sigma'_Z$ is invariant under $W$ and the induced
morphism from $\frak h/W$ to $Z+X$ is the Kostant section. In other words,
$\Sigma'_Z=\Sigma_Z$.
\end{lemma}
\begin{proof} Clearly $\Sigma'_Z(h)$ is conjugate to $h+X$, and hence by Part
(iii)  of Theorem~\ref{conj} maps to the same point of $\frak h/W$ as
$h$. Since there 
is a unique
such point in $Z+X$ with this property, namely
$\Sigma_Z (h)$, it 
follows that $\Sigma'_Z=\Sigma_Z$.
\end{proof}

We have thus proved (iv) of Theorem~\ref{conj}. Since the image of
the Kostant 
section lies in $\frak g_{\rm reg}$, it follows that the image of
$i+X$ has this 
property as well, proving (ii) of Theorem~\ref{conj}. This 
completes the proof of the theorem. \qed
\bigskip

\subsection{Consequences for invariant polynomials}

Let $\rho\colon G \to \Aut V$ be an irreducible representation and let 
$\rho_*\colon \frak g\to \operatorname{End} V$ be the induced representation of Lie
algebras. 
We can decompose $V$ into a direct sum of weight spaces $V_\nu$. Let $\mu$ be
the lowest 
weight of $\rho$, with corresponding weight space $V_\mu$, and let $\pi_\mu
\colon V\to V_\mu$ be the induced projection, i.e.\ $\pi_\mu$ is the identity on
$V_\mu$ and is zero on $V_\nu, \nu\neq \mu$. As in
\S2, let
$W_0$ be the stabilizer in $W$ of $\mu$.

The endomorphism $\rho_*(h_0)$ of $V$ is semisimple and its eigenvalues are
integral. Thus $V=\bigoplus _kV_k$, where $\rho_*(h_0)$ acts on $V_k$ with
eigenvalue $k$, and each $V_k$ is a sum of weight spaces. Clearly
$\rho_*(X)(V_k) \subseteq V_{k+2}$. The eigenspace $V_k$ corresponding
to the minimal value of $k$ is one dimensional, and in fact is exactly $V_\mu$. 

Let $S=\Sym^*\frak h^*$ be the affine coordinate ring of $\frak h$. The groups
$W$ and $W_0$ act on $S$, and we have the inclusions $S^W \subseteq
S^{W_0}\subseteq S$, corresponding to the finite surjective morphisms $\frak h
\to \frak h/W_0 \to \frak h/W$. Note  that the linear function $\mu$ is naturally
an element of
$S$ and in fact lies in $S^{W_0}$. Of course, $W$ acts on $\frak h$ as a group
generated by reflections. But,  by Lemma~\ref{generated},
$W_0$ also acts on $\frak h$ as a group generated by reflections.
Thus, by Chevalley's theorem  both $S^W$ and 
$S^{W_0}$ are polynomial algebras. We wish to compare these polynomial
algebras.

\begin{lemma}\label{faithflat} The ring $S$ is a faithfully flat extension of
$S^{W_0}$. Likewise,
$S^{W_0}$ is a faithfully flat extension of $S^W$. In particular, as an
$S^W$-module, $S^{W_0}$ is free of rank $\#(W/W_0)$.
\end{lemma}
\begin{proof} This is an immediate consequence of the fact that $S$, $S^{W_0}$,
and $S^W$ are regular and that the morphisms
in question are finite and surjective.
\end{proof}

For the rest of this section we fix a linear complement $Z\subseteq
\frak g$ to $\operatorname{Im}(\ad X)$, invariant under $h_0$, and we let $\Sigma
=\Sigma_Z$ and $\Psi$   be as given in Theorem~\ref{conj}.

\begin{lemma} 
Let $\lambda =\exp(-\Psi)\colon \frak h \to U_-$, so that
$\rho(\lambda)$ is a morphism from $\frak h$ to $\Aut V$. Then, for all $h\in
\frak h$, we have the following equality in $\operatorname{End} V$:
$$\rho(\lambda(h))\circ \rho_*(\Sigma(h)) = \rho_*(h+X)
\circ\rho(\lambda(h)).$$ 
\end{lemma} 

\begin{proof} This is immediate from Theorem~\ref{conj}.
\end{proof}

A morphism $\frak h\to \operatorname{End} V$ is the same thing as an element of
$$(\operatorname{End}_\Cee V)\otimes 
_\Cee S \cong \operatorname{End}_S (V\otimes _\Cee S).$$
Thus, after choosing a basis
of $V$, we 
can identify the  morphism $\rho(\lambda)$ from $\frak h$ to $V$ with an $n\times
n$  matrix with coefficients in $S$. The decomposition
of $V$ into the direct sum  of its weight spaces means that we can write
$\rho(\lambda)$ in block form as a sum of elements  
$\rho(\lambda)_{\nu_1, \nu_2}\in \Hom(V_{\nu_1}, V_{\nu_2})$, where $\nu_1,\nu_2$
are weights of $\rho$. It follows from Remark~\ref{homogeneous} that the entries
of  $\rho(\lambda)_{\nu_1, \nu_2}$ are homogeneous polynomials of the appropriate
degree. Similarly, the morphism $\ov \Sigma \colon \frak h/W \to \frak g$ induces
a morphism
$\rho_*\ov \Sigma\colon \frak h/W = \Spec  S^W  \to \operatorname{End} V$, and is thus
identified with an element of $\operatorname{End} V\otimes _\Cee S^W$, which we continue to
denote by $\rho_*\ov \Sigma$. The morphism $\rho_*\Sigma$ can likewise be viewed
as an element of $\operatorname{End} V\otimes _\Cee S$. Since $\Sigma$ is the pullback of $\ov
\Sigma$, it follows that $\rho_*\Sigma$ is the image of $\rho_*\ov \Sigma$ under
the inclusion $\operatorname{End} V\otimes _\Cee S^W\subseteq \operatorname{End} V\otimes _\Cee S$.

\begin{proposition} Define the $S$-submodule $M$ of $\Hom _S(V\otimes
_\Cee S,S)$ by
$$M =\{\phi\in \Hom _S(V\otimes _\Cee S,S): \phi \circ \rho_*(\Sigma)  =
\mu\cdot \phi\}.$$
Then $M$ is a free $S$-module of rank one, with a generator given by
$(\pi_\mu\otimes \Id)\circ \rho(\lambda)$ {\rm (}under any identification of
$V_\mu$ with $\Cee${\rm )}.
\end{proposition}

\begin{proof} Since $\rho_*(\Sigma)$ and $\rho_*(i+X)$ are conjugate in
$\Aut_S(V\otimes _\Cee S)$ (via $\rho(\lambda)$), it will suffice to prove the
same statement with $\rho_*(\Sigma)$ replaced by $\rho_*(i+X)$ and
$M$ replaced 
by 
$$M'=\{\phi\in \Hom _S(V\otimes _\Cee S,S): \phi \circ \rho_*(i+X)  =
\mu\cdot \phi\}.$$
In this case, $\phi\in M'$ if and only if $\phi\circ \rho(\lambda) \in M$.

Fix a basis $\{v_i\}$ of $V$ such that, for every $i, 1\leq i\leq n$, $v_i\in
V_{k_i}$. We may assume that, for $i\leq j$, $k_i\leq k_j$. In particular,
$v_1\in V_\mu$. In this basis, the element $\rho_*(i+X)$ is a lower triangular
matrix whose diagonal entries are the weights of $V$, viewed as
elements of $S$. 
Clearly, then, $M'$
is identified with the set of all $\phi\colon V\otimes _\Cee S \to S$ which
vanish on all weight spaces except for $V_\mu$. This space is  identified with
$V_\mu^*\otimes _\Cee S$, where $V_\mu^*\subseteq V^*$ is the
inclusion dual to the 
surjection $\pi_\mu\colon V \to V_\mu$. Thus, $M'$ is a free rank one
$S$-module 
with basis $\pi_\mu\otimes \Id$. It follows that $(\pi_\mu\otimes \Id)\circ
\rho(\lambda)$ is a basis for $M$.
\end{proof}

Let $W$ act on $S$ in the usual way and trivially on $V$. There are induced
actions of $W$ on $V\otimes _\Cee S$ and on $\Hom _S(V\otimes _\Cee S,S)$, via
$w\cdot \phi(v) = w\cdot \phi(w^{-1}v)$.

\begin{lemma}\label{W_0inv} The stabilizer $W_0$ of $\mu$ in $W$ acts trivially
on 
$(\pi_\mu\otimes \Id)\circ\rho(\lambda)$.
\end{lemma}

\begin{proof} By the previous proposition, $\tilde f=(\pi_\mu\otimes
\Id)\circ\rho(\lambda)$ is an eigenvector for the action of $\rho_*(\Sigma)$ on
$\Hom _S(V\otimes _\Cee S,S)$, with eigenvalue $\mu$, and the corresponding
eigenspace is free of rank one as an $S$-module. Since
$\rho_*(\Sigma)$ is
$W$-invariant, $w\cdot \tilde f$ is an eigenvector for $\rho_*(\Sigma)$ with
eigenvalue $w\cdot \mu$. In particular, if $w\in W_0$, then $w\cdot \tilde f=
c\tilde f$ for some $c\in S$. Since $\lambda\in G_-$, $\rho(\lambda)|V_\mu =
\Id$. Thus, 
$(\pi_\mu\otimes\Id)\circ\rho(\lambda)|V_\mu\otimes _\Cee S$ is the identity.
On the other hand, given $v\in V_\mu\otimes _\Cee S$ and $w\in W_0$,
clearly $w^{-1}v\in V_\mu\otimes _\Cee S$ as well, so that
$(w\cdot \tilde f)(v) = w\cdot w^{-1}v= v$. It follows that $w\cdot
\tilde f|V_\mu\otimes _\Cee S =\Id$, so that $c=1$. Hence $\tilde f$ is
$W_0$-invariant.
\end{proof}

\begin{corollary}\label{free} Let $M^{W_0}$ be the submodule of $M$ fixed by
$W_0$, so that 
$$M^{W_0} =\{\phi\in \Hom _{S^{W_0}}(V\otimes _\Cee S^{W_0}, S^{W_0}): \phi\circ
\rho_*(\Sigma) =\mu\cdot \phi\}.$$
Then $M^{W_0}$ is a free rank one $S^{W_0}$-module, generated by
$(\pi_\mu\otimes
\Id)\circ\rho(\lambda)$.
\end{corollary}

\begin{proof} We have a natural inclusion $M^{W_0} \otimes
_{S^{W_0}}S\to M$. The 
element $(\pi_\mu\otimes \Id)\circ\rho(\lambda)$ lies in $M^{W_0}$,
and so this 
inclusion is an isomorphism.
Thus, $M^{W_0}\otimes_{S^{W_0}}S$ is free of rank one.
 Since $S$ is faithfully flat over $S^{W_0}$, it
follows that  $M^{W_0}$ is a free rank one $S^{W_0}$-module and that 
$(\pi_\mu\otimes
\Id)\circ\rho(\lambda)$ is a generator.
\end{proof}

By Lemma~\ref{W_0inv}
the map $(\pi_\mu\otimes \Id)\circ\rho(\lambda)\colon V\otimes _\Cee
S\to S$
is  invariant under $W_0$ and hence induces a map
$$(\pi_\mu\otimes \Id)\circ\rho(\lambda)\colon V\otimes_\Cee
S^{W_0}\to S^{W_0}.$$
Let
$$f\colon  V\to S^{W_0}$$
be the restriction of $(\pi_\mu\otimes \Id)\circ\rho(\lambda)$ to
$V\otimes \{1\}\subset V\otimes_\Cee S^{W_0}$.
Of course, the map $f$ depends on $\lambda$,  which depends on the choice
of the slice $Z$.

\begin{corollary}\label{welldef}
Let $\frak m_0\subseteq S^W$ be the ideal generated by all
$W$-invariant 
polynomials whose constant term is zero. Then
the map
$$\ov f_0\colon V\to S^{W_0}/\frak m_0\cdot S^{W_0}$$
obtained by composing $f$ with the natural quotient mapping
$S^{W_0}\to
S^{W_0}/\frak m_0\cdot S^{W_0}$
is independent of the choice of slice $Z$.
\end{corollary}

\begin{proof}
If $Z'$ is another slice for $[X,\frak g]$, then by
Lemma~\ref{iso}
there is a map $\gamma\colon Z\to \frak u_-$ with $\gamma(0) =0$ such that
$$\Ad(\exp(\gamma(z)))(z+X)\in Z'+X$$
for all $z\in Z$. Let $\lambda=\exp(-\Psi) $ and $\lambda'=\exp(-\Psi')$, where
$\Psi,\Psi'$ are the maps given by Theorem~\ref{main} for the slices $Z$ and $Z'$,
respectively.
Then
$$\Ad(\exp(-\gamma(\Sigma_Z(h)-X)))\cdot
(\lambda')^{-1}\cdot\Ad(\exp(\gamma(\Sigma_Z(h)-X)))\cdot
\lambda
\colon  Z+X\to Z+X .$$ By uniqueness, it follows that 
$$\exp(\gamma(\Sigma_Z(h)-X))\cdot \lambda = \lambda'\cdot
\exp(\gamma(\Sigma_Z(h)-X)).$$ 
Now
$\exp(\gamma(\Sigma_Z(h)-X))$
is the exponential of a $W$-invariant map from $\frak h$ to $\frak u_-$, and thus
$\exp(\gamma(\Sigma_Z(h)-X)) \equiv 1\mod{\frak m_0}$. Hence,
$\lambda(h)\equiv
\lambda'(h)
\mod{\frak m_0\cdot S^{W_0}}$. 
\end{proof}

A Kostant section can be viewed as parametrizing a universal family of
regular, semi\-stable $G$-bundles on the cuspidal Weierstrass cubic $E$
whose support is contained in $E_{\rm reg}$.
We now apply Theorem~\ref{spectralcoverthm} and Theorem~\ref{spectralcoverthm2}
to  the 
family of vector bundles on $E$ obtained from this universal
family by a minuscule or quasiminuscule representation $\rho$.

\begin{theorem}\label{mainminu}
Let $\rho\colon G\to GL(V)$ be an irreducible representation with
lowest weight
$\mu$   and let
$$f=(\pi_\mu\otimes \Id)\circ \rho(\lambda)\colon V\to S^{W_0}$$
be the element given in Corollary~\ref{free}. Then by extension of
scalars $f$ induces a homomorphism
$\hat f\colon V\otimes_\Cee S^W\to S^{W_0}$
making the following diagram commute:
$$\begin{CD}\label{diagram}
V\otimes_\Cee S^W @>{\hat f}>> S^{W_0} \\
@V{\rho_*\ov\Sigma}VV  @VV{\mu\cdot}V \\
V\otimes_\Cee S^W @>{\hat f}>> S^{W_0}.
\end{CD}
$$
If  $\rho$ is minuscule, then $\hat f$ is an isomorphism. If  $\rho$ is
quasiminuscule, then $\hat f$ is surjective.
\end{theorem}

\begin{proof}
By Corollary~\ref{free} the  $S^{W_0}$-module map
$f\otimes \Id\colon V\otimes_\Cee S^{W_0}\to S^{W_0}$ makes the analogous diagram
commute, where $S^W$ is replaced by $S^{W_0}$ and $\rho_*\ov\Sigma$  by
$\rho_*\Sigma$, and the set of all such maps
$V\otimes_\Cee S^{W_0}\to S^{W_0}$ with this property is a free $S^{W_0}$-module
with $f\otimes \Id$ as a generator. Let $\hat f=f\otimes \Id|{V\otimes_\Cee S^W}$.
Since
$\rho_*(\Sigma)$ is
$S^W$-invariant we see that the diagram as defined in the theorem is
commutative.
It remains to show that $\hat f$ is surjective in the quasiminuscule case and an
isomorphism in the minuscule case.

Let $E$ be a cuspidal curve of arithmetic genus one.
We consider the $G$-bundle over $(\frak h/W)\times E$ which is trivial
on $(\frak h/W)\times \widetilde E$ and which is given by the Kostant
section
$\Sigma\colon\frak h/W\to \frak g_{\rm reg}\subset \frak g$ under the
correspondence of Theorem~\ref{triv}. Let $\mathcal{V}\to (\frak
h/W)\times E$ be 
the vector bundle induced from this $G$-bundle by $\rho$.
This vector bundle is also trivialized on $(\frak h/W)\times\widetilde
E$ and is given by the map $\rho_*\Sigma\colon\frak h/W\to \operatorname{End} (V)$.

Applying Theorem~\ref{spectralcoverthm} and Theorem~\ref{spectralcoverthm2} and
using the fact that every line bundle over $\frak h/W_0$ is trivial
 since $\frak h/W_0$ is isomorphic
to an affine space, we see  that there is a surjection from $\mathcal{V}$ to the
vector bundle
$$(\hat\nu\times \Id)_*(r\times\Id)^*\Poin,$$
where $\hat\nu\colon \frak h/W_0\to \frak h/W$ is the covering projection
and $r\colon \frak h/W_0\to \Cee\subseteq  E$ is the map induced
by the $W_0$-invariant  weight $\mu\colon E_{\rm reg}\otimes_\Zee
\Lambda\to E_{\rm reg}$. By Proposition~\ref{cuspPoin}, $\Poin\to E_{\rm
reg}\times E$ pulls back to the trivial line bundle on $E_{\rm
reg}\times\widetilde E$ and hence is given by a map $E_{\rm reg}\to  \operatorname{End}
(\Cee)=\Cee$. This map is the fixed identification $E_{\rm reg} \cong
\Cee$. Thus, 
when we use this identification to produce an isomorphism of the
coordinate ring of
$E_{\rm reg}$ with $\Cee[t]$, $\Poin$ becomes the
element $t\in 
\Cee[t]$. Hence, $(r\times\Id)^*\Poin$ is identified with the element
$r^*(t)$ which is the  element $\mu\in S^{W_0}$.
In this way we identify $(r\times \Id)^*\Poin$ with the
element
$\mu\in S^{W_0}$ viewed as an $S^{W_0}$-valued endomorphism of
$\Cee$.
Applying $(\hat\nu\times\Id)_*$ produces the $S^W$-module $S^{W_0}$
with the $S^W$-linear endomorphism given by multiplication by $\mu$.

Theorem~\ref{spectralcoverthm} now implies that there is a surjective $S^W$-linear
homomorphism
$f'\colon  V\otimes_\Cee S^W \to S^{W_0}$ making the following diagram commute:
$$\begin{CD}
V\otimes S^W @>{f'}>> S^{W_0} \\
@V{\rho_*\ov\Sigma}VV  @VV{\mu\cdot}V \\
V\otimes S^W @>{f'}>> S^{W_0}.
\end{CD}
$$
Its extension to a map $\hat f'\colon V\otimes_\Cee
S^{W_0}\to S^{W_0}$ is an 
element of the free module $M^{W_0}$ and hence $\hat f'=s(f\otimes \Id)$ for some
$s\in S^{W_0}$. Since $\hat f'$ is  surjective,
$s\in S^{W_0}$ is invertible.
Restricting to $V\otimes _\Cee S^W$, we see that $f'=s\hat f$,  and, since
$f'$ is surjective, so is
$\hat f$. Finally, if $\rho$ is minuscule, then $f'$ is an isomorphism and
hence $\hat f$ is an isomorphism as well.
\end{proof}

\begin{corollary}\label{cor} Suppose that $\rho$ is minuscule. 
Let $x\in \frak h/W$ and let $\frak m_x\subseteq S^W$ be the maximal
ideal of  $x$. Then the scheme-theoretic fiber
over $x$ in $\frak h/W_0$
has coordinate ring $S^{W_0}/\frak m_x\cdot S^{W_0}$. The map $\hat f$
induces an isomorphism $\ov f_x\colon V\to S^{W_0}/\frak m_x\cdot
S^{W_0}$. Under this isomorphism, multiplication by $\mu$ on the
right-hand side becomes the action of  $\rho_*\ov\Sigma(x)$. In particular, under
the isomorphism $\ov f_0\colon V\to S^{W_0}/\frak m_0\cdot
S^{W_0}$ given in Corollary~\ref{welldef},
multiplication by $\mu\in S^{W_0}$ on $S^{W_0}/\frak m_0\cdot
S^{W_0}$ corresponds to the
action 
$\rho_*(X)$ of the principal nilpotent element $X$ on $V$.\qed
\end{corollary}
 
In particular, the corollary describes the action of every regular element of
$\frak g$ on $V$ up to conjugation.

Using his deep results on perverse sheaves on the loop group of the Langlands
dual of $\frak g$, Ginzburg has proved the following generalization of 
Theorem~\ref{mainminu}:

\begin{theorem}\label{conj2} Let $\rho$ be an arbitrary irreducible
representation. The map $\hat f\colon V\otimes _\Cee S^W \to S^{W_0}$ is always
surjective.
\end{theorem}

It is natural to ask if there is a more elementary and direct proof of
Theorem~\ref{conj2}. The image of
$\hat f$ is the
$S^W$-submodule of $S^{W_0}$ generated by $(\pi_\mu\otimes \Id)\circ
\rho(\lambda)_{\nu, \mu}$, where $\rho(\lambda)_{\nu, \mu}\colon V_\nu \otimes
_\Cee S\to V_\mu \otimes _\Cee S$ is the linear map induced by $\rho(\lambda)$.
It follows the image of $\hat f$ is generated by homogeneous elements, when $S^W$
and $S^{W_0}$ are given their natural gradings as subrings of $S$. By the
homogeneous Nakayama lemma, to prove Theorem~\ref{conj2}, it would suffice to
prove that the induced homomorphism $\ov f_0\colon V\to S^{W_0}/\frak m_0\cdot
S^{W_0}$ is surjective.

\subsection{The Steinberg section}

Let  $\ov\Phi\colon H/W \to G_{\rm reg}$ be a 
section of Steinberg type. By this we mean that $\ov \Phi$ is any morphism from
$H/W$ to $G_{\rm reg}$, the set of regular elements of $G$, which is a section of
the adjoint quotient morphism
$G\to H/W$. We begin by proving a weaker analogue of Theorem~\ref{conj} for $\ov
\Phi$:

\begin{theorem} Let $\ov\Phi\colon H/W \to G_{\rm reg}$ be the Steinberg
section, and let $\Phi \colon H \to G_{\rm reg}$ be the composition $H\to H/W
\to G_{\rm reg}$. Then there exists a Borel subgroup $B$ of $G$ containing $H$,
with unipotent radical $U$, and a morphism
$\phi \colon H \to G$ such that, for all $h\in H$, $\phi(h)\Phi (h)
\phi(h)^{-1} \in B$, and moreover $\phi(h)\Phi (h)
\phi(h)^{-1} = hu(h)$ for some morphism $u\colon H\to U$.
\end{theorem}
\begin{proof} Suppose that we can show that there exists some Borel
subgroup $B$ and a morphism $\phi \colon H \to G$ such that, for all $h\in H$, 
$\phi(h)\Phi (h)
\phi(h)^{-1} \in B$. Then after a further conjugation we may assume that $B$
contains
$H$. In this case, we can write $\phi(h)\Phi (h)
\phi(h)^{-1} = h'u(h)$ for some morphism $u\colon H\to U$, where $h$ and $h'$ lie
in $H$ and have the same image in $H/W$. It follows that there is a $w\in W$ such
that, for an open dense set of $h\in H$, $h'= w(h)$. Hence we can assume that
$h'=w(h)$ for all $h$. After further conjugating $B$ by a representative for $w$
in the normalizer of $H$ in $G$, we can then assume that $\phi(h)\Phi (h)
\phi(h)^{-1} = hu(h)$ for some morphism $u\colon H\to U$ as desired.

To find the morphism $\phi$, we first claim that there exists a morphism
$\ov{\phi}$ from $H$ to the space $\mathcal{B}\cong G/B$ of Borel subgroups of
$G$, such that, for all $h\in H$, $h\in \ov{\phi}(h)$. To see this, fix a Borel
subgroup
$B$ and let
$I\subseteq G\times (G/B)$ be the incidence variety:
$$I =\{(g,xB): x^{-1}gx\in B\}= \{(g,xB): g\in xBx^{-1}\} .$$
Identifying the set of all Borel subgroups of $G$ with $G/B$, $I$ is the set of
pairs  consisting of an element $g$ of $G$ and a Borel subgroup
containing $g$. Clearly
$I$ is a closed subvariety of
$G\times (G/B)$. The morphism
$(g, xB)\mapsto (x, x^{-1}gx)$ defines an isomorphism from $I$ to 
$G\times _BB$, where
$B$ operates on
$G$ by right multiplication and on itself by conjugation, with inverse $(g,b)
\mapsto (gbg^{-1}, gB)$ (cf.\
\cite[\S4.3]{Slod}). Let
$I_{\rm reg}$ be the inverse image of $G_{\rm reg}$ under the projection
$\pi_1\colon I
\to G$. There is a morphism $\theta \colon I \to H$ defined by the
composition
$$I \cong G\times _BB \to G\times _BH \cong (G/B) \times H \to H,$$
where the morphism $G\times _BB \to G\times _BH$ is induced by the
homomorphism $B \to H$, and the morphism $(G/B) \times H \to H$ is projection
onto the second factor. It is straightforward to verify that the following
diagram is commutative:
$$\begin{CD}
I @>{\pi_1}>> G\\
@V{\theta}VV @VVV\\
H @>>> H/W.
\end{CD}$$
By a theorem of Grothendieck \cite[\S4.4]{Slod}, this diagram identifies
$I\to H$ with a simultaneous resolution of the the morphism $G\to H/W$, in
the terminology of \cite{Slod}. Since the adjoint quotient morphism
$G
\to H/W$ is smooth along
$G_{\rm reg}$, the above diagram identifies $I_{\rm reg}$ with the fiber product
$H\times _{H/W}G_{\rm reg}$.  Using the morphisms
$\ov\Phi
\colon H/W\to G$ and
$\pi_1\colon I\to G$, we can take the fiber product $(H/W)\times _GI$. By the
above remarks (all products below are fiber products), 
$$(H/W)\times _GI = (H/W)\times _{G_{\rm reg}}I_{\rm reg} \cong (H/W)\times
_{G_{\rm reg}}(H\times _{H/W}G_{\rm reg}) \cong H.$$
 The isomorphism $H \cong (H/W)\times _GI$ together with the  projection
$(H/W)\times _GI\to I$ identify an element
$h$ of $H$ with a pair $(x, B_h)$, where $x\in H/W$ is the image of $h$ and $B_h$
is a Borel subgroup containing $\ov\Phi(x)=\Phi(h)$. Here $B_h$ is the image of
$h$ under the morphism $\ov{\phi}\colon H \to G/B$, which is the
composition
$$H
\cong (H/W)\times _GI\to I
\to G/B,$$
and  we identify $G/B$ with the variety of all Borel subgroups of $G$.

Now suppose that we can lift the morphism $\ov{\phi}\colon H\to G/B$ to a morphism
$\phi\colon H \to G$. It then follows that, for all $h\in H$, $\phi(h)\Phi (h)
\phi(h)^{-1} \in B$ as claimed. To see that such a lift is possible, let
$\underline{B}$ be the sheaf of morphisms from $H$ to $B$, and similarly for
$\underline{G}$ and $\underline{G/B}$. Then the following is an exact sequence of
sheaves of sets in the \'etale topology:
$$\{1\} \to  \underline{B}\to  \underline{G} \to  \underline{G/B} \to\{1\}.$$
In particular, we have the following long exact sequence (of pointed sets):
$$H^0(H; \underline{G}) \to H^0(H; \underline{G/B}) \to H^1(H; \underline{B}).$$
On the other hand, the solvable group $B$ has a filtration by normal subgroups,
such that the successive quotients are either $\mathbb{G}_m$ or $\mathbb{G}_a$.
Since $H$ is a torus, all line bundles on it are trivial, and since it is affine,
all higher coherent sheaf cohomology vanishes. Thus $H^1(H; \underline{B})$ is
trivial. It follows that every morphism from $H$ to $G/B$ lifts to a morphism
from $H$ to $G$, as claimed. This concludes the proof.
\end{proof}

The method of proof above also proves a weak form of 
Theorem~\ref{conj}. However, it does not give the homogeneity statements.

Arguments very similar to those used in the proof of Theorem~\ref{mainminu} then
show the following:

\begin{theorem} Let $\mathcal{S}$ be the affine coordinate ring of $H$.
Let $\rho\colon G\to GL(V)$ be an irreducible representation with
lowest weight
$\mu$   and let
$$g=(\pi_\mu\otimes \Id)\circ \rho(\phi)^{-1}\colon V\to \mathcal{S}^{W_0}.$$
Then by
extension of scalars $g$ induces a homomorphism
$\hat g\colon V\otimes_\Cee \mathcal{S}^W\to \mathcal{S}^{W_0}$
making the following diagram commute:
$$\begin{CD}
V\otimes_\Cee \mathcal{S}^W @>{\hat g}>> \mathcal{S}^{W_0} \\
@V{\rho_*\ov\Phi}VV  @VV{\mu\cdot}V \\
V\otimes_\Cee \mathcal{S}^W @>{\hat g}>> \mathcal{S}^{W_0}.
\end{CD}
$$
If  $\rho$ is minuscule, then $\hat g$ is an isomorphism. If  $\rho$ is
quasiminuscule, then $\hat g$ is surjective. \qed
\end{theorem}

In case $\rho$ is minuscule, there is also a corollary  to the above theorem,
along the lines of Corollary~\ref{cor}, describing the action of
$\ov\Phi(x)$ on $V$ in terms of multiplication by $\mu$ on
$\mathcal{S}^{W_0}/\mathfrak{m}_x\cdot \mathcal{S}^{W}$, which we leave to the
reader to formulate.

Finally, it is natural to make the following conjecture, which does not seem to
follow from the Lie algebra analogue:

\begin{conjecture} In the above notation, the homomorphism $\hat g$ is always
surjective.
\end{conjecture}

\section{Examples}

\subsection*{The case of $A_{n-1}$}

Let $\frak h$ be  the Cartan subalgebra of $SL_n$, so that $\frak h\subset
\Cee ^n$ is identified with
$$\frak h=\left\{(x_1, \dots, x_n)|\sum _ix_i=0\right\}.$$
The
coordinate ring of $\frak h$ is
$S = \Cee[x_1, \dots, x_n]/(\sum _ix_i)$.
 Let
$\sigma_2, \dots, \sigma _n$ be the elementary symmetric functions in
$x_1,\ldots,x_n$ (note that $\sigma _1=0$ on $\frak h$).  Let $W =\frak S_n$
be the symmetric group on $n$ letters and let $W_0 =\frak S_{n-1}$, which is
embedded in $W$ as the stabilizer of $n$. Thus $W$ acts on $S$  and
$S^W =\Cee[\sigma_2, \dots, \sigma_n]$. Let
$\sigma_i(k)$
be the $i^{\rm th}$ elementary symmetric function in the variables $x_1,
\dots,x_k$. Thus for example $\sigma_i(n) =\sigma_i$. Note that $x_n
=-\sigma_1(n-1)$. More generally, we have the relation
\begin{equation*}\sigma_i(n-1)+x_n\sigma_{i-1}(n-1)=\sigma_i.\tag{$*$}
\end{equation*}
Clearly
\begin{align*}
S^{W_0} &=
\Cee[x_n,\sigma_1(n-1), \dots,
\sigma_{n-1}(n-1)]/(x_n+\sigma_1(n-1)=0)\\
&=\Cee[\sigma_1(n-1),\sigma_2(n-1),
\dots,
\sigma_{n-1}(n-1)].
\end{align*}
 By ($*$), it follows that $S^{W_0} = S^W[x_n]$. Since $x_n$ satisfies the monic
equation $p(x) = \sum _{i=0}^n(-1)^i\sigma_ix^{n-i} =0$, $S^{W_0} = S^W\cdot 1
\oplus \cdots \oplus S^W\cdot x_n^{n-1}$. Again by using ($*$), we can write
$x_n^i$ as $\sigma_i(n-1)$ plus terms involving $x_n^j$ for $j<i$ as well as
elements of $S^W$. Thus
$$S^{W_0} = S^W\cdot  1 \oplus S^W\cdot \sigma_1(n-1) \oplus \cdots \oplus
S^W\cdot \sigma_{n-1}(n-1)$$
as well, so that $1,\sigma_1(n-1), \dots, \sigma_{n-1}(n-1)$ is a basis for
$S^{W_0}$ as an $S^W$-module.

    Let
$V$ be the defining
$n$-dimensional representation of
$SL_n$. The vector space $V$ has the standard  basis  $e_1,\ldots,e_n$.
Define the
$n\times n$ matrix $\ov \Sigma$ by

$$\ov \Sigma = \begin{pmatrix}
0&1&0&0 &\cdots &0&0\\
0&0&1&0&\cdots  &0&0\\
0&0&0&1&\cdots  &0&0\\
\vdots&\vdots&\vdots&\vdots&\ddots &\vdots&\vdots\\
(-1)^{n-1}\sigma_n&(-1)^{n-2}\sigma_{n-1}&(-1)^{n-3}\sigma_{n-2}&
(-1)^{n-4}\sigma_{n-3}&\cdots &-\sigma_2&0
\end{pmatrix}.$$
As is well-known, $\ov\Sigma$ is a Kostant section for $\frak{sl}_n$.

For each $k,1\leq k\leq n$, let $m_i(k)$ be the sum of all monomials of degree
$i$ in the variables $x_1, \dots, x_k$. Define the column vector $\mathbf{c}_k$
by
$$\mathbf{c}_k = \begin{pmatrix}
0\\
0\\
\vdots\\
1\\
m_1(k)\\
\vdots\\
m_{n-k}(k)
\end{pmatrix},$$
and let
$A = (\mathbf{c}_1 \  \cdots \   \mathbf{c}_n)$.

It is not difficult to show that
$$A^{-1}\ov\Sigma A = \begin{pmatrix}
x_1&1&0&0 &\cdots &0&0\\
0&x_2&1&0&\cdots  &0&0\\
0&0&x_3&1&\cdots  &0&0\\
\vdots&\vdots&\vdots&\vdots&\ddots &\vdots&\vdots\\
0&0&0&0&\cdots&0&x_n
\end{pmatrix}$$
and that $A^{-1} =  $
$${\footnotesize
\begin{pmatrix}
1&0&0 &\cdots &0&0\\
-\sigma_1(1)&1&0&\cdots  &0&0\\
\sigma_2(2)&-\sigma_1(2)&1&\cdots  &0&0\\
\vdots&\vdots&\vdots&\ddots &\vdots&\vdots\\
(-1)^{n-1}\sigma_{n-1}(n-1)&(-1)^{n-2}
\sigma_{n-2}(n-1)&(-1)^{n-3}\sigma_{n-3}(n-1)&
\cdots&-\sigma_1(n-1)&1
\end{pmatrix}.}$$

Thus, $A^{-1}\in {\rm End}(V)\otimes S$ is the matrix realizing the
endomorphism $\rho(\lambda)$ of Theorem~\ref{main2}
for the standard $n$-dimensional representation  $\rho$ of $SL_n$.
The lowest weight of this representation is $x_n$ and hence
$\pi_\mu=x_n$.
Thus, $(\pi_\mu\otimes \Id)\circ \rho(\lambda)$ is the row matrix given by the
last row of $A^{-1}$.
Note that the last row of $A^{-1}$ is
invariant under $W_0=  \frak S_{n-1}$ and is in fact up to sign the basis for
$S^{W_0}$ over $S^W$ described above.

We now consider the minuscule representation
$\bigwedge^kV$. The lowest weight $\mu$ for this representation is
$e_n+\ldots +e_{n-k+1}$ and hence the map $(\pi_\mu\otimes \Id)\circ
\rho(\lambda)$ is the map given by the $k\times k$ minors of $A^{-1}$ which
involve the last $k$ rows.
The stabilizer of $\mu$ is $\frak S_{n-k}\times \frak S_k$,
where the first factor acts on the first $n-k$ variables and the
second on the last $k$ variables.
  It is clear from the explicit description of the
matrix for $A^{-1}$ that each of the last $k$ rows is invariant under the
group $\frak S_{n-k}$ and hence the determinants we are considering
are also invariant under this group.
 The theorem implies that these
determinants are also
invariant under $\frak S_k$ acting on the last $k$ variables.
In fact, this is  easy to check directly.

Let $\tau_j(k)$ be the $j^{\rm th}$ elementary symmetric function in the
last $k$ variables.
We have the relations
$$\sigma_\ell=\sum_{i=0}^\ell\sigma_i(n-k)\tau_{\ell-i}(k), \ell = 1, \dots,
n.$$ 
The ring $S^{W_0}$ is the quotient of the polynomial ring $\Cee[\sigma_i(n-k),
\tau_j(k)]$ modulo the above relations. The theorem implies
that the    $k\times k$-minors of the $k\times n$ matrix obtained from $A^{-1}$
by only considering the last $k$ rows are a basis for $S^{W_0}$ over $S^W$. 

To complete the picture for $SL_n$, we consider the adjoint representation. In
this case $\mu = e_n-e_1$ and the group $W_0$ is a copy of $\frak S_{n-2}$,
embedded in $\frak S_n$ as the stabilizer of $1$ and $n$. It is straightforward
to check in this case that the image of $V\otimes _\Cee S^W$ in $S^{W_0}$ is the
$S^W$-submodule of $S^{W_0}$ generated by $\sigma_{n-i}(n-1)x_1^{j-1}, 1\leq i,
j\leq n$, and that this submodule is all of $S^{W_0}$.

\subsection*{The case of $C_n$}

The unique minuscule representation of $Sp(2n, \Cee)$ is the defining
$2n$-dimensional representation.
Its weights   are $\pm x_i, 1\le i\le n$, with lowest weight $-x_1$.
The ring $S^W$ is a polynomial
algebra $\Cee[\alpha_1, \dots, \alpha_n]$, where $\alpha_i$ is the $i^{\rm th}$
elementary symmetric function  on  $x_1^2,\ldots,x_n^2$. If $W_0$ is the
stabilizer of the lowest weight
$-x_1$, then $S^{W_0}$ is a polynomial algebra with generators
$x_1,\beta_1,\ldots,\beta_{n-1}$ where
$\beta_i$ is the $i^{\rm th}$ elementary symmetric function on
$x_2^2,\ldots,x_n^2$. As in the case of $A_{n-1}$, there is a set of relations
$$x_1^2\beta_{i-1}+\beta_i =\alpha_i , i=1,\ldots,n.$$
Note that $x_1$ satisfies a monic polynomial over $S^W$ of degree $2n$.
It follows that $S^{W_0} = S^W[x_1] = S^W\cdot 1 \oplus S^W\cdot x_1 \oplus
\cdots \oplus S^W\cdot x_1^{2n-1}$. In particular, $-x_1$ acts on $S^{W_0}/\frak
m_0\cdot S^{W_0}$ with a Jordan block of length $2n$, corresponding to the fact
that the image of a principal nilpotent element of $\frak{sp}(2n)$ under the
standard representation is a principal nilpotent matrix. In terms of bundles, this
is equivalent to the statement that the bundle $I_{2n}$ has a nondegenerate
symplectic form.

\subsection*{The orthogonal representation of $Spin(2n)$}

Choose coordinates $x_1,\ldots, x_n$ on the Cartan subalgebra of
$\frak{spin}(2n)$ so that the coroot lattice is the even integral lattice.
The orthogonal representation of $Spin(2n)$ has  weights $\pm
x_1,\ldots,\pm x_n$, and  $-x_1$ is the lowest weight.  Let $W_0$ be the subgroup
of the Weyl group $W$ stabilizing
$-x_1$. It is the Weyl group for $Spin(2n-2)$ in the variables  $x_2,\ldots,
x_n$. As before, we let $\alpha_i$ be the $i^{\rm th}$
elementary symmetric function  on  $x_1^2,\ldots,x_n^2$ and $\beta_i$ the $i^{\rm
th}$ elementary symmetric function on
$x_2^2,\ldots,x_n^2$. The ring $S^W$ is the
ring generated by 
$\alpha_1, \dots, \alpha_n$
 and the Pfaffian
${\rm Pfaff}_n=x_1\cdots x_n$ subject to the one relation
${\rm Pfaff}_n^2=\alpha_n$. Thus $S^W$ is a polynomial 
algebra on
$\alpha_i, 1\le i <n$, and
${\rm Pfaff}_n$. The algebra $S^{W_0}$ is the algebra generated by
$x_1$, $\beta_i,  1\le i\leq n-1$, and ${\rm Pfaff}_{n-1}$,
subject to the relation ${\rm Pfaff}_{n-1}^2=\beta_{n-1}$, 
and thus
$S^{W_0}$ is a polynomial algebra generated by
$x_1,\beta_i, 1\le i<n-1,$ and ${\rm Pfaff}_{n-1}$.
The  ring
$S^{W_0}/\frak m_0\cdot S^{W_0}$
is the quotient of $S^{W_0}$  by the relations
\begin{align*}
x_1^2\beta_i+\beta_{i+1}&=0;\\
x_1{\rm Pfaff}_{n-1}&=0.
\end{align*}
Hence $S^{W_0}/\frak m_0\cdot S^{W_0}$ is generated by $x_1$ and ${\rm
Pfaff}_{n-1}$ subject to the relations
\begin{align*}
&x_1^{2n-2}={\rm Pfaff}_{n-1}^2;\\
&x_1{\rm Pfaff}_{n-1}=0.
\end{align*}
It follows that $$x_1^{2n-1}=x_1{\rm Pfaff}_{n-1}^2=0.$$
A $\Cee$-basis for $S^{W_0}/\frak m_0\cdot S^{W_0}$ is
$1,x_1,\ldots,x_1^{2n-2},{\rm Pfaff}_{n-1}$. As a module
over $\Cee[x_1]$, $S^{W_0}/\frak m_0\cdot S^{W_0}\cong
\Cee[x_1]/(x_1^{2n-2})\oplus
\Cee[x_1]/(x_1)$. This reflects the fact, which is easy to check directly, that
a  principal nilpotent element of
$\frak{spin}(2n)$ acts on
$V$ with two Jordan blocks, one of dimension $2n-1$ and one
of dimension $1$. This implies that
the orthogonal bundle of rank $2n$ which is regular as a $Spin(2n)$-bundle and
S-equivalent to the trivial bundle
is $I_{2n-1}\oplus \mathcal{O}_E$. In particular, it is not
regular as a vector bundle.

\subsection*{The case of $E_6$}

Direct computation shows that the action of the principal nilpotent
element $X$ in the Lie algebra of type  of $E_6$ on the  $27$-dimensional
minuscule representation
$V$   has three Jordan blocks of dimensions $17$, $9$, and $1$.
The bottom of the Jordan block of rank $17$ is the lowest weight space
$V_\mu$. Thus, Theorem~\ref{main2} implies that multiplication by
$\mu$ in $S^{W_0}/\frak m_0\cdot S^{W_0}$ is nilpotent of order $17$.
 One can show that $S^{W_0} = S^W[\mu, a, b]$, where
$\mu, a,b$ are homogeneous in degrees $1$, $4$, $8$ respectively. Note that  
$S^W$ is a graded polynomial algebra with generators in degrees $2, 5, 6, 8, 9,
12$. On the other hand, since $W_0$ is the Weyl group of $D_5$, and its action on
$\frak h$ decomposes as a copy of the standard representation of $W(D_5)$ plus a
trivial factor, it follows that
$S^{W_0}$ has generators in degrees
$1, 2,4, 5, 6 ,8$. Given these degrees and the fact that $S^{W_0}$ has rank $27$
over $S^W$, it is easy to check directly that $S^{W_0}$ is generated as an
$S^W$-algebra by elements in degrees $1,4,8$.

\section{Quasiminuscule representations}

Let $\rho$ be a quasiminuscule, nonminuscule representation, and let
$\mathcal{V}$ be the corresponding vector bundle over $\Mod \times E$. In this
section, we describe the extension class corresponding to the extension
\begin{equation*}
0 \to \pi_1^*\pi_1{}_*\mathcal{V} \to \mathcal{V} \to
\overline{\mathcal{V}} 
\to 0 \tag{$*$}
\end{equation*}
of Theorem~\ref{spectralcoverthm2}, in the case where $G$ is simply laced and
hence $\rho$ is the adjoint representation. It seems likely that similar
geometric methods can also describe the quasiminuscule, nonminuscule
representation in the non-simply laced case. We shall just write out the case of
$\frak g$, i.e.\ where $E$ is cuspidal, although minor modifications handle the
nodal and smooth cases. Recall the notation of \S2.1: we have the finite cover
$\hat
\nu\colon \hat{T}_0 \to \Mod$, as well as the normalization map $\varphi\colon
\hat{T}_0 \to T_0 \subseteq \Mod \times E_{\rm reg}$, where $T_0$ is the image of
$(\hat \nu, r) \colon \hat{T}_0\to \Mod\times E_{\rm reg}$. Fixing  a 
coordinate $t$ for $E$ at $p_0$, the pullback
$r^*t$ in the coordinate ring of $\hat{T}_0$, i.e.\ as an element of $S^{W_0}$, is
just the weight
$\mu$. In the above notation, $\ov{\mathcal{V}} = (\hat \nu\times
\Id)_*(r\times \Id)^*\Poin$.

\subsection{The extension group} We begin by identifying the module of extensions
for exact sequences whose first and third terms are above.

\begin{proposition}\label{identifyext} Let $\mathcal{W}$ be a vector bundle
over $\Mod$ and let $M$ be the projective and hence free
$S^W$-module corresponding to
$\mathcal{W}$. Then $\Ext^1_{\Mod\times
E}(\overline{\mathcal{V}},\pi_1^*\mathcal{W})$, which is a module over $S^W$, is
isomorphic to the $S^W$-module $(S^{W_0}/\mu S^{W_0})\otimes _{S^W}M$.
\end{proposition}
\begin{proof} By Corollary~\ref{extcor}, the module of such extensions is
isomorphic to 
$(Q\spcheck/tQ\spcheck)\otimes _{S^W}M$, where $Q\spcheck$ is the Fourier-Mukai
transform of $(\overline{\mathcal{V}})\spcheck$. Thus it suffices to show that
$Q\spcheck \cong S^{W_0}$. By Theorem~\ref{spectralcoverthm2},
$\overline{\mathcal{V}}$ is the direct image under $\hat \nu\times \Id$ of the
pullback of a line bundle on
$\hat{T}_0$, necessarily trivial, tensored with $(r\times
\Id)^*\Poin$. Applying duality for the finite flat morphism $\hat \nu\times \Id$,
and using the fact that the relative dualizing sheaf $\omega_{\hat{T}_0/\Mod}$ is
a line bundle on $\hat{T}_0$, and hence is trivial as well, it follows that 
$(\overline{\mathcal{V}})\spcheck$ is the direct image under $\hat \nu\times \Id$
of $(r\times \Id)^*\Poin^{-1}$. Note that reflection in the root $\mu$ defines an
element of
$W$ of order two normalizing $W_0$, and hence an involution $\tau$ of $\hat{T}_0$
such that $\hat\nu\circ\tau =\hat\nu$ and  $r\circ \tau = -\Id\circ r$.
Thus  
\begin{align*}(\hat \nu\times \Id)_*(r\times \Id)^*\Poin^{-1} &= (\hat \nu\times
\Id)_*(\tau\times \Id)^*(r\times \Id)^*\Poin^{-1} \\
&= (\hat \nu\times
\Id)_*(r\times \Id)^*(-\Id\times \Id)^*\Poin^{-1}.
\end{align*}
On the other hand, from the symmetric definition of $\Poin$, it follows that
$$(-\Id\times \Id)^*\Poin^{-1} \cong \Poin.$$ Thus
$(\overline{\mathcal{V}})\spcheck\cong (\hat \nu\times
\Id)_*(r\times \Id)^*\Poin \cong\overline{\mathcal{V}}$, and in particular
$Q\spcheck \cong S^{W_0}$.
\end{proof}

 Let $D = T_0 \cap (\Mod\times \{p_0\})$ and let $\hat D$ be the
pullback of $D$ to $\hat{T}_0$. Since the affine coordinate ring of $\hat D$ is
equal to $S^{W_0}/\mu S^{W_0}$, we can rephrase the proposition as follows:

\begin{corollary} In the above notation,  the $S^W$-module $\Ext^1_{\Mod\times
E}(\overline{\mathcal{V}},\pi_1^*\mathcal{W})$ is isomorphic to the $S^W$-module
of global sections of $\hat\nu^*\mathcal{W}|\hat D$. \qed
\end{corollary}

Clearly, $\hat D = r^{-1}(0)$, and hence $\hat D$ is isomorphic to $\Mod_0$, in
the terminology of \S2.1. In particular, $\hat D$ is smooth. It follows from
Lemma~\ref{singofmap} that the induced morphism
$\hat D \to D$ is birational. We shall say more about the structure of $\hat D \to
D$ below.

\subsection{Identification of the extension class} 

Assume now that  $\rho$ is the adjoint representation. The nonzero weights of
$\rho$ are the roots of $\frak g$, and
$\mathcal{V} =\ad
\Xi$.  First, we identify  the vector bundle $\mathcal{W} = \pi_1{}_*\mathcal{V}$:

\begin{proposition} The vector bundle $\pi_1{}_*\mathcal{V}$ is isomorphic to
$\Omega^1_{\Mod}$.
\end{proposition}
\begin{proof} The bundle $\Xi$ is a universal $G$-bundle over $\Mod \times E$. It
follows that the Kodaira-Spencer homomorphism $T_{\Mod} \to R^1\pi_1{}_*\ad \Xi$
is an isomorphism. (Of course, this is just a restatement of the result of
Kostant \cite{Kostant} that the adjoint quotient morphism is smooth at a regular
element.) Using relative duality, the isomorphism $(\ad \Xi)\spcheck \cong \ad
\Xi$ given by the Killing form, and the fact that the relative dualizing sheaf
$\omega_{\Mod \times E/\Mod}$ is trivial, it follows that $(R^1\pi_1{}_*\ad
\Xi)\spcheck \cong R^0\pi_1{}_*\ad \Xi$. Thus 
$\Omega^1_{\Mod} \cong \pi_1{}_*\ad \Xi = \pi_1{}_*\mathcal{V}$, as claimed.
\end{proof}

There is another description of  $\pi_1{}_*\mathcal{V}\cong \Omega^1_{\Mod}$
as follows:

\begin{lemma}
The sheaf $\pi_1{}_*\mathcal{V}\cong \Omega^1_{\Mod}$
corresponds to the  $S^W$-module $(S\otimes _\Cee \frak h)^W$, where $W$ acts on
both factors of the tensor product by the natural actions.
\end{lemma}
\begin{proof} The sheaf $\Omega^1_{\Mod}$ is the coherent sheaf associated to
the $S^W$-module $\Omega^1_{S^W/\Cee}$. By  a theorem of Solomon (see for example
\cite[p.\ 135, ex.\ 3]{Bour}),
$\Omega^1_{S^W/\Cee} = (\Omega^1_{S/\Cee})^W$. But $\Omega^1_{S/\Cee}\cong
S\otimes _\Cee \frak h^*$, and the invariant bilinear form defines a
$W$-equivariant isomorphism from $\frak h^*$ to $\frak h$. Thus
$\Omega^1_{S^W/\Cee} \cong (S\otimes _\Cee \frak h)^W$.
\end{proof}

The extension class corresponding to $\mathcal{V}$ is thus a global section  of
$\hat\nu^*\Omega^1_{\Mod}|\hat D$. Our goal now will be to describe this
section. Viewing
$D$ as a divisor on $\Mod \times
\{p_0\} \cong \Mod$, there is the conormal morphism
$I_D/I_D^2\to  \Omega^1_{\Mod}|D$. The image of this morphism is the conormal
sheaf. At a smooth point of $D$, this subsheaf is a  line subbundle which is the
dual of the normal bundle to $D$. Since
$\Mod$ is an affine space,
$\Pic
\Mod$ is trivial, and hence there is an $\bar s \in S^W$, unique up to a nonzero
scalar, such that $I_D = (\bar s)$. In fact, we can take $\bar s = \prod_{w\in
W/W_0}w\mu$. Thus
$d\bar s$ defines a section of
$\Omega^1_{\Mod}|D$. There is then the induced section
$\hat\nu^*d\bar s$ of $\hat\nu^*\Omega^1_{\Mod}|\hat D$. Since $\hat D \cong
\Mod_0$, which is again an affine space, there exists an indivisible section $s$
of
$\hat\nu^*\Omega^1_{\Mod}|\hat D$, unique up to a constant, and an element $a\in
\Gamma (\hat D,
\scrO_{\hat D})\cong S^{W_0}/\mu S^{W_0}$ such that $\hat\nu^*d\bar s = as$.
Concretely, choosing a basis of the trivial vector bundle $\Omega^1_{\Mod}$ and
hence of 
$\hat\nu^*\Omega^1_{\Mod}|\hat D\cong \scrO_{\hat D}^r$, we can write $s= (f_1,
\dots, f_r)$, where the $f_i$ have no common factor, and $as = (af_1, \dots,
af_r)$ for functions
$f_i$ on
$\hat D$.

\begin{theorem}\label{classext} The above section
$s$ defines the extension class corresponding to $\mathcal{V}$. Moreover, $s$ can
be completed to a basis for $\hat\nu^*\Omega^1_{\Mod}|\hat D$. In other words,
the sheaf $(\hat\nu^*\Omega^1_{\Mod}|\hat D)/\scrO_{\hat D}\cdot s$ is a locally
free sheaf of $\scrO_{\hat D}$-modules.
\end{theorem}
\begin{proof} We begin by verifying that, at a generic point $x$ of $D$, the
corresponding extension class lies in the image of the conormal line $N_{D/\Mod,
x}^{-1}$ in $\Omega^1_{\Mod, x}$. For such an 
$x\in D$, we may assume that
$D$ and $T_0$ are smooth at $x$. Let $\xi$ be the $G$-bundle over $E$
corresponding to $x\in \Mod$. By the proof of Theorem~\ref{spectralcoverthm} and
\ref{spectralcoverthm2}, $\xi$ is described as follows: let $H'$ be the subtorus 
of $H$ given by $\Ker
\mu$ and let $SL_2$ be the subgroup of $G$ corresponding to the Lie algebra
spanned by $\frak g^\mu$,
$\frak g^{-\mu}$, and $\mu\spcheck$. Then there is the natural embedding 
$SL_2\times _{\Zee/2\Zee}H' \to G$. Let $\xi_1$ be the $SL_2$-bundle 
corresponding to
$I_2$. Then there exists  a generic $H'$-bundle $\xi_2$ such that, if
$\xi'$ is the bundle on 
$SL_2\times _{\Zee/2\Zee}H'$ induced by $\xi_1\boxtimes \xi_2$ on $SL_2\times 
H'$, then  $\xi$ is the induced $G$-bundle. Note that $\ad_{SL_2}\xi_1 =I_3$ and
that $\ad_{H'}\xi_2 =\scrO_E^{r-1}$. There is the induced injection 
$\ad_{SL_2}\xi_1\oplus \ad_{H'}\xi_2 \to \ad_G\xi$. The arguments of
Theorem~\ref{spectralcoverthm2} imply that the corresponding homomorphism on
$H^1$ is an isomorphism. In particular,
$$T_{\mathcal{M}, x} \cong H^1(E; \ad_G \xi) = H^1(E; \ad_{SL_2}\xi_1) \oplus 
H^1(E;
\ad_{H'}\xi_2).$$
Clearly, $H^1(E; \ad_{H'}\xi_2)$ is  the image of the tangent space $T_{D,x}$ to
$D$ at  $x$, and hence the normal bundle sequence is split at $x$. Another way to
give the same splitting on the normal bundle sequence  is as follows. The
Killing form induces a quadratic form on
$H^1(E;  
\ad \xi)$, which is easily seen degenerate with radical equal to $H^1(E;
\ad_{SL_2}\xi_1)$, and this splits the sequence.
Dualizing, there is a direct sum decomposition 
$$\Omega^1_{\mathcal{M}, x} \cong H^0(E; \ad \xi) = H^0(E; \ad_{SL_2}\xi_1) 
\oplus H^0(E;\ad_{H'}\xi_2).$$
Under the identification of $\Omega^1_{D, x}$ with $H^0(E;\ad_{H'}\xi_2)$, the
above is a splitting of the conormal sequence such that the natural morphism
$\Omega^1_{\mathcal{M}, x} \to
\Omega^1_{D, x}$ is the projection to $H^0(E;\ad_{H'}\xi_2)$. 

The extension $(*)$ restricted to the slice $\{x\}\times E$ gives an extension
$$0\to H^0(E; \ad \xi) \otimes _\Cee\scrO_E \to \ad \xi \to \ov V \to 0,$$
and hence defines an extension class $\varepsilon\in H^1(E; \ov V{}\spcheck)
\otimes H^0(E;
\ad
\xi)$. By construction, the part of $\ov V$ supported at $p_0\in E$ is isomorphic
to
$I_2$ and  so 
$H^1(E; \ov V{}\spcheck) \cong \Cee$. Thus up to a scalar we can view
$\varepsilon$ as an element  of $H^0(E; \ad \xi)$. We claim that, if $V'$ is
defined via the commutative diagram
$$\begin{CD}
0 @>>> H^0(E; \ad \xi) \otimes _\Cee\scrO_E @>>> \ad \xi @>>> \ov{V} @>>> 0\\
@. @VVV @VVV @V{=}VV @.\\
0 @>>> H^0(E; \ad_{H'} \xi_2) \otimes _\Cee\scrO_E @>>> V' @>>> \ov{V} @>>> 0
\end{CD}$$
then the corresponding extension $V'$ is split. There are many ways to see this. 
For example, it follows from the explicit knowledge of $\ad \xi$ given in the
proof of Theorem~\ref{spectralcoverthm2}. Another way is to use that fact that
orthogonal projection from $\frak g$ to $\frak h'$ is equivariant with respect
to the action of the subgroup $SL_2\times _{\Zee/2\Zee}H'$. Thus there is an
induced homomorphism of vector bundles from $\ad \xi$ to
$H^0(E;
\ad_{H'} \xi_2) \otimes _\Cee\scrO_E$. In turn, this homomorphism induces a
surjection from
$V'$ to $H^0(E; 
\ad_{H'}
\xi_2) \otimes _\Cee\scrO_E$ which splits the extension. In any case, we see that
$\varepsilon\in H^0(E; \ad \xi) \cong \Omega^1_{\mathcal{M}, x}$
lies in the kernel of the projection $\Omega^1_{\mathcal{M},
x} \to \Omega^1_{D, x}$, and hence $\varepsilon$ lies in the image of the
conormal line
$N_{D/\Mod, x}^{-1}$. Thus, at a generic point of $D$, the extension class
$\varepsilon$ is equal to a nonzero multiple of $\hat{\nu}^*d\bar s$ at $x$ and
hence is a nonzero multiple of $s$ at $x$.

Suppose that  
$s'$ is the section of $\hat \nu^*\Omega^1_{\Mod}|\hat D$ corresponding to
the extension $(*)$. The above shows that there exists a regular function $c$ on
$\hat D$ such that
$s' = cs$.  Fixing a
basis for
$\Omega^1_{\Mod}$, we can write
$s' = cs = (cf_1, \dots, cf_r)$.  To complete the
proof of Theorem~\ref{classext}, it suffices to show the following:   for all
$y\in
\hat D$, there exists an
$i$ such that $cf_i\notin \frak m_y$. This will show that $c$ is a nowhere
vanishing function on $\hat D$ and hence a constant, and that $s$ is a nowhere
vanishing section of  $\hat \nu^*\Omega^1_{\Mod}|\hat D$. To see this, let
$x\in \Mod$ lie  under $y$, let
$\xi$ be the corresponding
$G$-bundle, and let $0\to \scrO_E^r \to V \to \ov V \to 0$ be the corresponding
extension, where we use the given basis of $\Omega^1_{\Mod}$ to trivialize 
the subbundle $H^0(E;V)\otimes \scrO_E$ of $V$. The fact about $V$ that we shall
need is that $h^0(E; V)=r$ and hence that the coboundary map $\delta \colon
H^0(E; \ov V) \to H^1(E; \scrO_E^r)$ is injective. Write $\delta = (\delta_1,
\dots, \delta_r)$, where   $\delta_i \in \Ext^1(\ov V, \scrO_E) \cong \Hom(H^0(E;
\ov V), H^1(E;\scrO_E))$. By Proposition~\ref{identifyext}, we have identified
$Q=Q(V)$ with
$\widetilde{R} = S^{W_0}/\frak m_xS^{W_0}$ and
$\Ext^1_{\scrO_E}(\ov V, \scrO_E^r)$ with $(\widetilde{R}/\mu\widetilde{R})^r$.
Under this identification, the $i^{\textrm{th}}$ component $\delta_i$ of $\delta$
corresponds to the $i^{\textrm{th}}$ component of $s'$, namely $cf_i$. Thus, if
$cf_i\in \frak m_y$ for every $i$, then $cf_i$ lies in the subspace $\frak
m_y/\mu\widetilde{R}$ of $\widetilde{R}/\mu\widetilde{R}$, which is a proper
subspace since $y\in \hat D$, i.e.\ $\frak m_y$ is not the unit ideal.  This
would imply that  all of the $\delta_i$ lie in a proper subspace of $\Hom(H^0(E;
\ov V), H^1(E;\scrO_E))$, and hence that they have a common kernel. This
contradicts the injectivity of $\delta$. Thus $cf_i \notin \frak m_y$ for some
$i$, proving that both $c$ and $f_i$ are units at $y$. 
\end{proof}

We can say more about the function $a$ such that $\hat\nu^*d\bar s = as$ in the
discussion preceding  Theorem~\ref{classext}. To determine $a$ up to a nonzero
constant, it suffices to describe the divisor of zeroes  of
$a$ in $\hat D$. We begin by describing the divisor
$D$ in more detail. Let
$y\in
\hat D$, so that by definition $y$ is the image of a point $e\in E_{\rm
reg}\otimes
\Lambda$ such that $\mu(e)=0$. If $\pm\mu$ are the only roots  which vanish at
$e$, then by Lemma~\ref{singofmap}, $\hat{T}_0 \to T_0$ is a local isomorphism
near
$y$, and hence $\hat D\to D$ is also an isomorphism. Thus we may assume that
there exists a root $\mu'\neq \pm \mu$ such that $\mu'(e) =0$ as well. The
possibilities are as follows (we leave the calculations to the reader):

\begin{lemma}\label{hatD} Let $e\in E_{\rm reg}\otimes
\Lambda$ lie over $y\in \hat D$. Suppose that there exists a root $\mu'\neq \pm
\mu$ such that $\mu'(e) =0$, and that $e$ is a generic point of $D(\mu)\cap
D(\mu')$.
\begin{enumerate}
\item[\rm (i)] If $\mu$ and $\mu'$ are orthogonal, then $\varphi^{-1}(\varphi(y))
=\{y, y'\}$ consists of two distinct points, both $T_0$ and $D$ have two local
branches at
$\varphi(y)$, meeting transversally, and $\varphi\colon \hat{T}_0
\to T_0$ is a local diffeomorphism at $y$ onto one of the two branches of $T_0$
passing through $\varphi(y)$. There exist
local coordinates $x_1, \dots, x_{r-2}$, $s_1$, $t_2$ for $\hat{T}_0$ at $y$
and local coordinates $x_1, \dots, x_{r-2}$, $t_1$, $t_2$ for $\Mod$ at $\hat
\nu(y)$  such that
$$\hat \nu(x, s_1, t_2) = (x, t_1, t_2),$$
where $t_1=s_1^2$.  A local form for the equation of $D$ is $t_1t_2 =0$. Finally,
$\mu = r^*t = s_1$.
\item[\rm (ii)] If $\mu$ and $\mu'$ are not orthogonal, then they span a root
system of type $A_2$. The preimage $\varphi^{-1}(\varphi(y)) =\{y\}$. There exist
local coordinates $x_1, \dots, x_{r-2}$, $s_1$, $s_2$ for $\hat{T}_0$ at $y$ and 
local coordinates $x_1, \dots, x_{r-2}$, $\sigma_2$, $\sigma_3$ for $\Mod$ at
$\hat
\nu(y)$  such
that
$$\hat\nu(x, s_1, s_2) = (x, \sigma_2, \sigma_3),$$
where $\sigma_2 = -(s_1^2+s_2^2 + s_1s_2)$ and $\sigma_3 = -(s_1^2s_2 + s_1
s_2^2)$. A local form for the equation of $D$ is $4\sigma_2^3 + 27\sigma_3^2 =0$.
Finally, $\mu =r^*t = s_1-s_2$.
\qed
\end{enumerate}
\end{lemma}

\begin{corollary} Let $\hat{D}_{1,1}$ the the hypersurface in $\hat D$ which is
the  image of the set of $e \in E_{\rm reg}\otimes
\Lambda$ such that there exists a root $\mu'\neq \pm
\mu$, orthogonal to $\mu$, with $\mu'(e) =0$. Let $\hat{D}_2$ the the hypersurface
in $\hat D$ which is the  image of the set of $e \in E_{\rm reg}\otimes
\Lambda$ such that there exists a root $\mu'\neq \pm
\mu$,  not orthogonal to $\mu$ with $\mu'(e) =0$. Then the function $a$ vanishes
to order
$1$ along every component of $\hat{D}_{1,1}$ and to order $3$ along every
component of $\hat{D}_2$.
\qed
\end{corollary}

\section{Further conjectures}

In this final section, we speculate on how some of the previous results and
conjectures might be generalized to an arbitrary irreducible representation
$\rho$.  For simplicity, we shall only discuss the case of $\frak g$, where these
questions might be accessible to the methods of Ginzburg. The general idea is to
find a filtration of the corresponding vector bundle $\mathcal{V}$, or
equivalently of the module $V\otimes _\Cee S^W$, so that the the associated
graded module is free and action of
$\rho_*\overline{\Sigma}$ on the associated graded can be described explicitly.
For simplicity, we shall shift our point of view and work with highest weights
instead of lowest weights in what follows.

\subsection{Notation}

Given an irreducible representation $\rho\colon G \to GL(V)$,
let
$\lambda$ be a weight of $\rho$ and let $m_\lambda$ be the multiplicity of
$\lambda$, so that
$V=\bigoplus _\lambda V_\lambda$, where $H$ acts on $V_\lambda$ by the character
corresponding to $\lambda$, and $\dim V_\lambda = m_\lambda$. Recall that, if
$\lambda$ is a weight, then so is $w\lambda$ for every $w\in W$. In particular,
each $W$-orbit of weights contains a unique dominant weight. We define a partial
ordering on the set of all dominant weights as follows: If $\lambda_1$ and
$\lambda_2$ are two dominant weights, then   $\lambda_1
\geq
\lambda_2$ if the difference $\lambda_1-\lambda_2$ is a rational linear
combination of simple roots with nonnegative coefficients, and 
$\lambda_1>\lambda_2$ to mean $\lambda_1\geq \lambda_2$ and $\lambda_1\neq
\lambda_2$.

\subsection{Statement of the conjecture}

For each dominant weight
$\lambda$ of
$\rho$, define the operator 
$P_\lambda \in \End (V\otimes_\Cee S^W)$ by
$$P_\lambda =\prod _{\lambda'\in W\cdot \lambda}(\rho_*\ov\Sigma - \lambda'\cdot
\Id).$$ Note that $P_\lambda$ does in fact lie in $\End (V\otimes_\Cee S^W)$,
although the individual  factors only lie in $\End (V\otimes _\Cee S)$, and that
$P_\lambda$ commutes with
$\rho_*\ov 
\Sigma$. Moreover, $\Ker P_\lambda$ is a saturated submodule of $V\otimes_\Cee
S^W$   which corresponds to the component of the spectral cover whose
reduction is given by the Weyl orbit of $\lambda$. It seems natural, however,
that the ordering of the weights must somehow be taken into account. To do so,
make the following construction: For
$\lambda$ dominant, define
$$Q_\lambda =\prod _{\lambda' \leq \lambda}P_\lambda,$$ 
where the product is over all dominant weights $\lambda'\leq \lambda$, define and
$$Q_\lambda^0 =\prod _{\lambda' < \lambda}P_\lambda.$$ 
Then $F_\lambda = \Ker Q_\lambda\subseteq V\otimes_\Cee
S^W$ is invariant under $\rho_*\ov 
\Sigma$, as is $F_\lambda^0 = \Ker Q_\lambda^0$, and $F^0_\lambda \subseteq 
F_\lambda$. Note that, since the dominant weights of $\rho$ might only be
partially ordered, $\{F_\lambda\}$ need not be a filtration of $V\otimes_\Cee
S^W$.

\begin{conjecture} Let $\lambda$ be a dominant weight. Let $W_\lambda$ be  the
stabilizer in
$W$ of $\lambda$ and $m_\lambda$  the multiplicity of $\lambda$ in $\rho$. Then,
in the above notation, there is an isomorphism of $S^W$-modules
$$\hat h \colon F_\lambda/F_\lambda^0 \cong (S^{W_\lambda})^{m_\lambda},$$  
and $\hat h\circ \rho_*\ov \Sigma = \lambda \hat h$.
\end{conjecture}

In case $\lambda$ is the highest weight of $\rho$, the conjecture is Ginzburg's
result. At the other extreme, if $\lambda =0$, the conjecture is contained in the
following:

\begin{proposition} $\Ker \rho_*\ov \Sigma$ is a free submodule of $V\otimes_\Cee
S^W$ of rank $m_0$.
\end{proposition}
\begin{proof} If $0$ is not a weight of $\rho$, then $\Ker \rho_*\ov
\Sigma(x)=0$ for each $x$ such that $\ov \Sigma(x)$ is semisimple. Thus $\Ker
\rho_*\ov \Sigma=0$, and hence it is free of rank $m_0=0$. Thus we may assume that
$0$ is a weight of $\rho$. In this case, the proof uses the following result of
Kostant
\cite[(5.1.1)]{Kostant}: Suppose that $m_0$, the multiplicity of the weight $0$ in
$V$, is positive. Let
$X$ be a regular element of $\frak g$, let  $\frak g^X$ denote the centralizer of
$X$, and define 
$$V^{\frak g^X}=\{v\in V: \rho_*(Y)(v) = 0 \mbox{ for all $Y\in \frak g^X$}\}.$$
Then $\dim V^{\frak g^X}=m_0$, and in particular it is independent of the choice
of $X$. 

For each $x\in \frak h/W$, we just write $\frak g^x$ for  $\frak g^{\ov \Sigma
(x)}$. Then $\bigcup_{x\in \frak h/W}\frak g^x$ is an $r$-dimensional subbundle
$\frak z_{\ov \Sigma}$ of the trivial bundle $(\frak h/W)\times \frak g$. Since
$S^W$ is a polynomial algebra, the vector bundle $\frak z_{\ov \Sigma}$
corresponds to a free $S^W$-submodule of $\frak g\otimes_\Cee S^W$, with basis
$X_1,
\dots, X_r$ (although the main point is rather that the submodule is projective). 
Let $M = \bigcap_{i=1}^r\Ker \rho_*(X_i) \subseteq V\otimes _\Cee S^W$. Since
$\dim V^{\frak g^x}$ is independent of $x$, it follows that $M$ is a projective,
and hence free, submodule of $V\otimes_\Cee S^W$, clearly saturated. As
$\ov\Sigma(x)
\in \frak g^x$ for every $x$, $M\subseteq \Ker \rho_*\ov \Sigma$. Finally, if $x$
is the image of a regular semisimple element, then it is easy to see that
$V^{\frak g^x} =\Ker\ov \Sigma(x)$, and since $M$ is saturated, we must have
equality everywhere. Thus $\Ker \rho_*\ov \Sigma =M$, and in particular $\Ker
\rho_*\ov \Sigma$ is free.
\end{proof}

\end{document}